\documentclass[12pt,twoside]{amsart}

\usepackage[dvipsnames]{xcolor}
\usepackage{hyperref}
\hypersetup{colorlinks=true,
allcolors=Blue}
\usepackage{amsthm, amsmath, amscd, amssymb,centernot,txfonts}
\usepackage{tikz-cd}
\usepackage{lipsum}
\usepackage{amsfonts}
\usepackage{mathrsfs}
\usepackage[all]{xy}
\usepackage[left=2cm,top=2.5cm,bottom=3cm,right=3cm]{geometry}
\usepackage{dirtytalk}

\theoremstyle{plain}
\numberwithin{equation}{section}

\newcommand{\myeq}{\overset{\mathrm{duality}}{=\joinrel=}}
\newcommand*{\QEDA}{\hfill\ensuremath{\blacksquare}}
\newcommand*{\QEDB}{\hfill\ensuremath{\square}}
\begin{document}

\title[Projective Normality And Normal Presentation on Certain Varieties]{On The Projective Normality And Normal Presentation On Higher Dimensional Varieties with Nef Canonical Bundle}

\author[Mukherjee]{Jayan Mukherjee}
\address{Department of Mathematics, University of Kansas, Lawrence, KS 66045}
\email{j899m889@ku.edu}

\author[Raychaudhury]{Debaditya Raychaudhury}
\address{Department of Mathematics, University of Kansas, Lawrence, KS 66045}
\email{debaditya@ku.edu}


\maketitle

\begin{abstract}
In this article we prove new results on projective normality and normal presentation of adjunction bundle associated to an ample and globally generated line bundle on higher dimensional smooth projective varieties with nef canonical bundle. As one of the consequences of the main theorem, we give bounds on very ampleness and projective normality of pluricanonical linear systems on varieties of general type in dimensions three, four and five. These improve known such results.

\end{abstract}

\section*{\textbf{Introduction}}
Equations defining the embedding of a projective variety in a projective space is a topic of great interest. The study of projective normality and normal presentation dates back to the time of Italian geometers. Castelnuovo first showed that a line bundle of degree greater than $2g$ on a curve of genus $g$ has a normal homogeneous coordinate ring and if the degree is greater than $2g+1$ then the ideal of the curve is generated by quadrics. Fujita, St. Donat and Mumford, among many others, rediscovered these results years later. Mumford and his school of mathematicians carried on the study of these properties on an abelian variety of aribitrary dimension. In the early 80s, Green and Lazarsfeld showed that the results of these nature are special cases of a general $N_p$ property (see \cite{G1}, \cite{G2} and \cite{G3}) for curves.\par 
We start with the definition of projective normality, normal presentation and the property $N_p$.

\noindent{\bf Definition 0.1.} \textit{Let $L$ be a very ample line bundle on a variety $X$. Let the following be the minimal graded free resolution of the coordinate ring $R$ of the embedding of $X$ induced by the complete linear system} $|L|$\[
\begin{tikzcd}
  0 \arrow{r}{} & F_n \arrow{r}{\phi_n} & F_{n-1} \arrow{r}{\phi_{n-1}} & ... \arrow{r}{\phi_1} & F_0\arrow{r}{\phi_0} & R \arrow{r}{} & 0.
\end{tikzcd}\]
\textit{Let $\mathscr{I}_X$ be the ideal sheaf of the embedding.
\begin{itemize}
    \item[(1)] $L$ satisfies the property $N_0$ (or embeds $X$ as a projectively normal variety) if $R$ is normal.
    \item[(2)] $L$ satisfies the property $N_1$ (or is normally presented) if in addition $\mathscr{I}_X$ is generated by quadrics.
    \item[(3)] $L$ satisfies the property $N_p$ if in addition to satisfying the property $N_1$, the resolution
is linear from the second step until the $p$-th step.
\end{itemize}}\par

Mark Green proved that a line bundle of degree $\geq 2g+1+p$ on a smooth curve of genus $g$ satisfies the property $N_p$. One of the most interesting questions on surfaces concerning the $N_p$ property that has motivated lot of work is \textit{Mukai's Conjecture : For an ample line bundle $A$ on a smooth projective surface $S$, $K_S+lA$ satisfies the $N_p$ property if $l\geq p+4$ ($K_S$ is the canonical bundle on $S$)}. This can be thought of as an analogue of Green's result on curves for surfaces.\par
Mukai's conjecture has not yet been proved even for $p=0$. Note that by Reider's result we have that $K_S+lA$ is very ample if $l \geq 4$ (see \cite{Re}). We state some of the results obtained on specific varieties towards this direction below.\par

\textit{Elliptic Ruled Surfaces:} Y. Homma proved it for the case $p=0$ for elliptic ruled surface (see \cite{H1} and \cite{H2}). The case $p=1$ for elliptic ruled surfaces was proved by Gallego and Purnaprajna. In fact, they showed that the numerical classes of normally presented divisors on an elliptic ruled surface forms a convex set and as a particular case recovered Mukai's conjecture for $p=0,1$ and yield weaker bounds for higher syzygies (see \cite{GP4}).\par
    
\textit{Ruled Varieties:} Butler proves that in characteristic $0$, if $E$ is a rank $n$ vector bundle on a smooth projective curve $C$ with genus $g \leq 1$ then $K_X+lA$ is projectively normal for $l\geq 2n+1$ and satisfies the property $N_p$ for $l \geq 2n(p+1)$ where $X=\mathbb{P}(E)$  (see \cite{Bu}).\par
     
\textit{Surfaces with Kodaira Dimension zero:} Gallego and Purnaprajna proved Mukai's conjecture on these surfaces for $p=0,1$ lowering the bound by one in the latter case (see \cite{GP6}).\par
    
\textit{Abelian Varieties:} On abelian varieties Koizumi's theorem states that the $lA$ is projectively normal for $l \geq 3$ and $A$ ample (see \cite{Koi}). Kempf further proved that  $lA$ is normally presented for $l \geq 4$ and $A$ ample (see \cite{Ke}). The above results on abelian varieties were generalized by Pareschi where he showed that $lA$ satisfies the property $N_p$ for $l\geq p+3$ (see \cite{Pa}).\par
     
\textit{Hyperelliptic Varieties:} Chintapalli and Iyer proved that $lA$ satisfies $N_p$ property for $l\geq p+3$ where $A$ is any ample line bundle on a hyperelliptic variety (see \cite{CI}).\par
    
\textit{Surfaces of General Type:} B.P. Purnaprajna proved that under mild hypothesis on an ample and globally generated line bundle $A$, $K+lA$ is projectively normal and normally presented for $l\geq 2$ where $K$ is the canonical line bundle. He also obtained precise results on higher syzygies (See \cite{BP}). $N_p$ property of the adjoint bundles associated to an ample and globally generated line bundle on surfaces of general type was also studied in \cite{BH} where Banagere and Hanumanthu proved several interesting results in this direction.\par
     
\textit{Toric Varieties:} Hering, Schenck and Smith proved in \cite{HSS} that for an ample line bundle $A$ on an $n$ dimensional toric variety, $lA$ satisfies $N_p$ property for $l\geq n+p-1$. \par 

Ein and Lazarsfeld proved that for a very ample line  bundle $L$ on a smooth projective variety $X$,  $K_X + (n+1+p)A$ satisfies the property $N_p$.  (see \cite{EL1}).\par
Another very interesting and related conjecture is the conjecture by Fujita. The precise statement is the following:\par
\textit{Fujita's Conjecture: On a smooth projective variety of dimension $n$, $K_X+(n+1)A$ is globally generated and $K_X+(n+2)A$ is very ample where $A$ is an arbitrary ample line bundle}.
\par
 Fujita's conjecture has been proved for surfaces by Reider (cf. \cite{Re}) using Bogomolov's instability theorem (see \cite{Bo}) on rank two vector bundles. Fujita's freeness conjecture has been proved by Ein and Lazarsfeld (see \cite{EL}) for $n=3$, by Kawamata (see \cite{Ka}) for $n=3,4$ and by Fei Ye and Zhixian Zhu (see \cite{Zhu}) for $n=5$. \par
Mukai's conjecture can be generalized as follows: \textit{For a smooth projective variety of dimension $n$ and an ample line bundle $A$, $K_X+lA$ satisfies the property $N_p$ for $l \geq n+p+2$}. Progress in this direction with $A$ just ample seems to be out of reach at this moment. A natural question to ask is what happens to the above conjecture if $A$ is taken to be ample and base point free instead. It is a standard argument that if $A$ is taken to be ample and base point free then Fujita's conjecture follows in its full generality by using induction and using known results for curves.
Syzygies of adjunction bundles with $A$ ample and base point free was studied in quite some details on surfaces in a series of papers written by Gallego and Purnaprajna (see \cite{GP1}-\cite{GP6}).\par 
In this paper we prove new results on the properties $N_0$ and $N_1$ of the adjunction bundle $K+lB$ with $B$ ample and base point free on arbitrary dimensional smooth projective varieties with nef canonical bundle by imposing mild conditions on the line bundle $B$ apart from the ones mentioned above. These are analogues for results known for surfaces. Our main result regarding projective normality on a variety $X$ with canonical divisor $K_X$ is the following:\par
 \noindent{\bf Theorem.} (See Theorem \hyperref[2.3]{2.3}) \textit{Let $X$ be a smooth projective variety of dimension $n$, $n\geq 3$. Let $B$ be an ample and base point free line bundle on $X$. We further assume:\\
\indent (a) $K_X$ is nef, $K_X+B$ is base point free.\\
\indent (b) $h^0(B)\geq n+2$.\\
\indent (c) $h^0(K_X+B)\geq h^0(K)+n+1$.\\
\indent (d) $B-K_X$ is nef and effective.\\
Then $K_X+lB$ is very ample and it embeds $X$ as a projectively normal variety for all $l\geq n$.} \par
Note that, in general $h^0(B)\geq n+1$ and $h^0(K_X+B)\geq h^0(K_X)+n$ (See Remark \hyperref[2.2.01]{2.2.1}) and hence the conditions are not as strong since $B$ is ample and base point free. 
Using the mildness of the conditions we impose on $B$ we come up with the following corollary: \par
\noindent{\bf Corollary.} (See Corollary \hyperref[2.4]{2.4}) \textit{Let $X$ be a variety of dimension $n\geq 3$ with $p_g\geq 2$. Let $B$ be an ample, globally generated line bundle on $X$. Assume $K_X$ is nef and $B-K_X$ is a nef, non-zero, effective divisor. Further assume that $B+K_X$ is globally generated. If either $H^1(B)=0$ or $H^{n-1}(\mathscr{O}_X)=0$ then $K_X+nB$ will be very ample and it will embed $X$ as a projectively normal variety.}\par

\textit{Sharpness of our conditions:} To discuss the sharpness of our conditions we produce two sets of examples. \par 
In Example \hyperref[2.5]{2.5} we produce examples of smooth projective varieties in all dimensions satisfying all conditions of Theorem \hyperref[2.3]{2.3} excepting $h^0(B) \geq n+2$ and show that $K+nB$ is not projectively normal, where $n$ is the dimension of the variety, thereby emphasizing the sharpness of the condition in the theorem. \par

In Example \hyperref[2.6]{2.6} we produce examples of smooth projective varieties in all dimensions that satisfy all conditions in Corollary \hyperref[2.4]{2.4} excepting the fact that $B-K$ is nef, non-zero and effective and show that $K+nB$ is not projectively normal, where $n$ is the dimension of the variety, thereby  showing that the condition is essential.\par

Our result regarding normal presentation is the following:\par

 \noindent{\bf Theorem.} (See Theorem \hyperref[3.4]{3.4} and \hyperref[3.5]{3.5}) \textit{Let $X$ be a smooth projective variety of dimension $n$ with nef canonical bundle, $n\geq 3$. Let $B$ be an ample and base point free line bundle on $X$ with $h^0(B)\geq n+2$. We further assume:\\
\indent (a) $K+B$ is base point free. In addition, if $X$ is irregular then for any line bundle $B'\equiv B$, $B'$ and \indent\indent$\:K+B'$ are base point free.\\
\indent (b) $h^0(K+B)\geq h^0(K)+n+1$. In addition, if $X$ is irregular then for any line bundle $K'\equiv K$, \indent\indent$\:h^0(K+B)\geq h^0(K')+n+1$.\\
\indent (c) $(n-2)B-(n-1)K$ is nef and non-zero effective divisor.\\
Then $K+lB$ will satisfy the property $N_1$ for $l\geq n$}.\par

Once we have these theorems, we can start looking for results using only an ample bundle if we know what multiple of that bundle is globally generated. Here solution to \textit{Fujita's freeness Conjecture} comes to play an important role. \par

The geometry of pluricanonical maps is of great importance in projective algebraic geometry. It was extensively studied by Bomberi, Catanese, Ciliberto, Kodaira (see \cite{Bom}, \cite{CC1}, \cite{CC2}, \cite{Kod}). Ciliberto showed that for minimal surfaces of general type $nK$ is projectively normal for $n \geq 5$ (see \cite{Ci}). B.P Purnaprajna produced very precise and optimal bounds for normal generation and normal presentation and higher syzygies of pluricanonical series on surfaces of general type with ample canonical bundle (see \cite{BP}). 

In this paper we obtain effective results on projective normality and normal presentation of pluricanonical series on smooth threefolds, fourfolds and fivefolds with ample canonical bundle. The following corollary is the summary of Corollaries \hyperref[4.3]{4.3}, \hyperref[4.4]{4.4}, \hyperref[4.5]{4.5}, \hyperref[4.6]{4.6}, and \hyperref[4.7]{4.7}:\par

 \noindent{\bf Corollary.}  \textit{Let $X$ be a smooth projective variety of dimension $n$ with ample canonical bundle $K$.\\
\indent (i) If $n=3$, then $lK$ is very ample and embeds $X$ as a projectively normal variety for $l \geq 12$ and \indent\indent normally presented for $l \geq 13$.\\
\indent (ii) If $n=4$, then $lK$ is very ample and embeds $X$ as a projectively normal variety for $l \geq 24$ and \indent\indent normally presented for $l \geq 25$.\\
\indent (iii) If $n=5$ and $p_g(X)\geq 1$, then $lK$ is very ample and embeds $X$ as a projectively normal variety \indent\indent for $l \geq 35$ and normally presented for $l \geq 36$}.
\par
As far as we know this corollary gives new bounds on very ampleness, projective normality and normal presentation of pluricanonical systems on threefolds, fourfolds and fivefolds. \par
The standard arguments using Castelnuovo-Mumford regularity yield weaker results.\par 
For example, for a smooth projective threefold with ample canonical bundle $K$, Castelnuovo-Mumford regularity shows that $nK$ satisfies projective normality and normal presentation for $n\geq 14$ and $n\geq 16$ respectively. So we need more subtle methods. We build on the methods of \cite{BP} and use newer ideas, one such is to use the Skoda complex.   \par
In the last section we generalize our results to projective varieties with Du-Bois singularites and hence derive some effective results on projective normality and normal presentation of pluricanonical series on projective threefolds with $\mathbb{Q}$-factorial terminal Gorenstein singularities or with canonical Gorenstein singularities. \par
\pagebreak
\textbf{Acknowledgements}. We are extremely grateful to our advisor Prof. B.P. Purnaprajna for introducing us to this subject, teaching us the key concepts and guiding us throughout this work. We also thank the referee for  suggestions and corrections that substantially helped in the improvement of the exposition. 

\section{Preliminaries and notations}\label{prelims}
Throughout this paper, we will always work on a projective variety $X$ over an algebraically closed field of characteristic zero. $K$ or $K_X$ will denote its canonical bundle. We will use the multiplicative and the additive notation of line bundles interchangeably. Thus, for a line bundle $L$, $L^{\otimes r}$ and $rL$ are the same. We have used the notation $L^{-r}$ for $(L^{*})^{\otimes r}$. We will use $L^{r}$ to denote the intersection product. The sign \say{$\equiv$} will be used for numerical equivalence.   \par 
Let $X$ be a smooth, projective variety and let $L$ be a globally generated line bundle on $X$. We define the bundle $M_L$ as follows\phantomsection\label{0},
\begin{equation}
    0 \longrightarrow M_L \longrightarrow  H^0(L)\otimes \mathscr{O}_X \longrightarrow L \longrightarrow 0. \tag{$\ast$}
\end{equation}

If L is an ample and globally generated line bundle on $X$ one has the
following characterization of the property $N_p$.\par
\noindent\textit{{\bf Theorem 1.1.}\phantomsection
\label{1.1} Let $L$ be an ample, globally generated line bundle on $X$. If the group $H^1(\bigwedge^{p'+1}M_L\otimes L^{\otimes k})$ vanishes for all $0\leq p'\leq p$ and for all $k\geq 1$, then L satisfies the property $N_p$. If in addition $H^1(L^{\otimes r})=0$ for all $r\geq 1$, then the above vanishing is a necessary and sufficient condition for $L$ to satisfy $N_p$.}\par
Since we are working over a field with characteristic zero, $\bigwedge^{p'+1}M_L$ is a direct summand of $M_L^{\otimes p'+1}$ (see \cite{EL1}, Lemma 1.6). Consequently, to show that a line bundle $L$ satisfies the property $N_p$, we will show that $H^1(M_L^{\otimes p'+1}\otimes L^{\otimes k})=0$ for all $0\leq p'\leq p$ and for all $k\geq 1$. Notice that $L$ being projectively normal automatically implies that $L$ is very ample. We refer to \cite{Mu} for an overview of these circle of ideas.\par
The following observation has been used often in the works of Gallego and Purnaprajna (see for instance \cite{GP6}).\par
\noindent\textit{{\bf Observation 1.2.}\phantomsection
\label{1.2} Let $E$ and $L_1$, $L_2$,..., $L_r$ be coherent sheaves on a variety $X$. Consider the map  $H^0(E)\otimes H^0(L_1\otimes L_2\otimes...\otimes L_r)\xrightarrow[]{\psi}H^0(E\otimes L_1\otimes...\otimes L_r)$ and the following maps
$$H^0(E)\otimes H^0(L_1)\xrightarrow[]{\alpha_1}H^0(E\otimes L_1),$$
$$H^0(E\otimes L_1)\otimes H^0(L_2)\xrightarrow[]{\alpha_2}H^0(E\otimes L_1\otimes L_2),$$
$$...$$
$$ H^0(E\otimes L_1\otimes ... \otimes L_{r-1})\otimes H^0(L_r)\xrightarrow[]{\alpha_r}H^0(E\otimes L_1\otimes...\otimes L_r).$$
If $\alpha_1$, $\alpha_2$,..., $\alpha_r$ are surjective then $\psi$ is also surjective.}\par

The following from \cite{GP2} relates the surjectivity of a multiplication map on a variety to the surjectivity of its restriction to a divisor.\par

\noindent\textit{{\bf Lemma 1.3.}\phantomsection
\label{1.3} Let $X$ be a regular variety (i.e. $H^1(\mathscr{O}_X)=0$). Let $E$ be a vector bundle and let $D$ be a divisor such that $L=\mathscr{O}_X(D)$ is globally generated and $H^1(E\otimes L^*)=0$. If the multiplication map $H^0(E\vert_D)\otimes H^0(L\vert_D)\rightarrow H^0((E\otimes L)\vert_D)$ surjects then $H^0(E)\otimes H^0(L)\rightarrow H^0(E\otimes L)$ also surjects.}\par

The proposition below is a result from \cite{Bu}. Here $\mu$ denotes the slope of a vector bundle.\par

\noindent\textit{{\bf Proposition 1.4.}\phantomsection
\label{1.4} Let $E$ and $F$ be semistable vector bundles over a curve $C$ of genus $g$ such that E is generated by its global sections. If\\
\indent (1) $\mu(F)> 2g$, and\\
\indent (2) $\mu(F)>2g+rank(E)(2g-\mu(E))-2h^1(E)$.\\
Then the multiplication map $H^0(E)\otimes H^0(F)\rightarrow H^0(E\otimes F)$ surjects.}\par

The following lemma from \cite{GP2} is an useful tool for showing normal presentation.\par
\noindent\textit{{\bf Lemma 1.5.}\phantomsection\label{1.5} Let $X$ be a projective variety, let $r$ be a non-negative integer and let $F$ be a base-point-free line bundle on $X$. Let $Q$ be an effective line bundle on $X$ and let $q$ be a reduced and irreducible member of $|Q|$. Let $R$ be a line bundle and $G$ a sheaf on $X$ such that\\
\indent (1) $H^1(F\otimes Q^*)=0$\\
\indent (2) $H^0(M_{F\otimes  \mathscr{O}_q}^{\otimes i}\otimes R\otimes \mathscr{O}_q)\otimes H^0(G)\rightarrow H^0(M_{F\otimes \mathscr{O}_q}^{\otimes i}\otimes R\otimes G\otimes\mathscr{O}_q)$ is surjective for all $0\leq i\leq r$.\\
Then for all $0\leq i'\leq r$ and for all $0\leq k\leq i'$, the following map
\begin{align*}
    H^0(M_F^{\otimes k}\otimes M_{F\otimes \mathscr{O}_q}^{\otimes i'-k}\otimes R \otimes \mathscr{O}_q)\otimes H^0(G)\longrightarrow H^0(M_F^{\otimes k}\otimes M_{F\otimes \mathscr{O}_q}^{\otimes i'-k}\otimes R \otimes G\otimes \mathscr{O}_q)
\end{align*}
is surjective.}\par 

The lemma below, a generalization of the base point-free pencil trick, is due to Green (c.f. \cite{G2}, Theorem (4.e.1)):\par

\noindent\textit{{\bf Lemma 1.6.}\phantomsection
\label{1.6} Let $C$ be a smooth, irreducible curve. Let $L$ and $M$ be line bundles on $C$. Let $W$ be a base point free linear subsystem of $H^0(C,L)$. Then the multiplication map $W\otimes H^0(M)\rightarrow H^0(L\otimes M)$ is surjective if $h^1(M\otimes L^{-1})\leq dim(W)-2$.}\par

The following lemma called the Castelnuvo-Mumford lemma (see \cite{Mu}) will be used frequently in this article.\par

\noindent\textit{{\bf Lemma 1.7.}\phantomsection
\label{1.7} Let $L$ be a base point free line bundle on a variety $X$ and let $\mathscr{F}$ be a coherent sheaf on $X$. If $H^i(\mathscr{F}\otimes L^{-i})=0$ for all $i\geq 1$ then the multiplication map $H^0(\mathscr{F}\otimes L^{\otimes i})\otimes H^0(L)\rightarrow H^0(\mathscr{F}\otimes L^{\otimes i+1})$ surjects for all $i\geq 0$.}\par

If the variety is not regular, we will not be able to use Lemma \hyperref[1.3]{1.3} to show the surjection of a multiplication map. To overcome the problem, we have to use the Skoda complex which is defined below. We will use it often to show the projective normality and the normal presentation on an arbitrary variety.\par

\noindent{\bf Definition 1.8.} \phantomsection
\label{1.8}\textit{Let $X$ be a smooth projective variety of dimension $n\geq 2$. Let $B$ be a globally generated and ample line bundle on $X$.\\
\indent (1) Take $n-1$ general sections $s_1,...s_{n-1}$ of $H^0(B)$ so the intersection of the divisor of zeroes \indent\indent $B_i=(s_i)_0$ is a nonsingular projective curve $C$, that is $C=B_1\cap ...\cap B_{n-1}$.\\
\indent (2) Let $\mathscr{I}$ be the ideal sheaf of $C$ and let $W=span\{s_1,...,s_{n-1}\}\subseteq H^0(B)$ be the subspace \indent\indent spanned by $s_i$. Note that $W\subseteq H^0(B\otimes \mathscr{I})$.
For $i\geq 1$, define the Skoda complex ${\bf I}_i$ as}
\[
\begin{tikzcd}[row sep=large, column sep=2ex]
 0 \arrow{r}{} & \bigwedge\limits^{n-1}W\otimes B^{-(n-1)}\otimes \mathscr{I}^{i-(n-1)} \arrow{r}{} & \dots\arrow{r}{} & W\otimes B^{-1}\otimes \mathscr{I}^{i-1}\arrow{r}{} & \mathscr{I}^i \arrow{r}{} & 0
\end{tikzcd}\]
\indent\indent\textit{where $\mathscr{I}^k$ stands for $\mathscr{I}^{\otimes k}$, we have used the convention that $\mathscr{I}^k=\mathscr{O}_X$ for $k\leq 0$.}\par 

In this article we have only used ${\bf I}_1$ which is the following,
\[
\begin{tikzcd}
  0 \arrow{r}{} & \bigwedge\limits^{n-1}W\otimes B^{-(n-1)} \arrow{r}{} & \dots\arrow{r}{} & \bigwedge\limits^{2}W\otimes B^{-2} \arrow{r}{} & W\otimes B^{-1}\arrow{r}{} & \mathscr{I} \arrow{r}{} & 0
\end{tikzcd}\] \\
and it is just the Koszul resolution of $\mathscr{I}$.

Even though our main theorems deal with the adjoint bundle associated to an ample and globally generated line bundle, in Section \ref{KX} we deduce some results on three, four and five folds that deal with the pluricanonical series when the canonical bundle is just an ample line bundle. In order to make this transition, we need Fujita's freeness conjecture on three, four and five folds or a slightly stronger version of it (see \cite{EL} and \cite{Ka}). In particular, we need the following results.\par

\noindent\textit{{\bf Theorem 1.9.} (See \cite{Ka}, Theorem 3.1)\phantomsection\label{1.9} Let $X$ be a normal projective variety of dimension 3, $L$ an ample Cartier divisor, and $x_0 \in X$ a smooth point. Assume that there are positive numbers $\sigma_p$ for $p = 1, 2, 3$ which satisfy the following conditions:\\
\indent (1) $\sqrt[\leftroot{-1}\uproot{1}p]{L^p\cdot W}\geq \sigma_p$ for any subvariety $W$ of dimension $p$ which contains $x_0$.\\
\indent (2) $\sigma_1\geq 3$, $\sigma_2\geq 3$ and $\sigma_3>3$.\\
Then $|K_X+L|$ is free at $x_0$.}\par

\noindent\textit{{\bf Corollary 1.10.} (See \cite{Ka}, Corollary 3.2)\phantomsection\label{1.10} Let $X$ be a smooth projective variety of dimension 3, and $H$ an ample divisor. Then $|K_X +mH|$ is base point free if $m \geq 4$. Moreover, if $H^3\geq 2$, then $|K_X + 3H|$ is also base point free.}\par

\noindent\textit{{\bf Theorem 1.11.} (See \cite{Ka}, Corollary 4.2)\phantomsection\label{1.11} Let $X$ be a smooth projective variety of dimension 4, and $H$ an ample divisor. Then $|K_X + mH|$ is base point free if $m\geq 5$.}\par
The remark after the following result from \cite{Mi} will be used in Section \ref{KX}.\par
\noindent\textit{{\bf Theorem 1.12.} Let $k$ be an algebraically closed field of characteristic 0 and $X$ a normal projective $\mathbb{Q}$-Gorenstein variety of dimension $n\geq 2$ with singular locus of codimension $\geq 3$. Assume that the canonical divisor $K_X\in Pic(X)\otimes \mathbb{Q}$ is nef. Let $\rho:Y\rightarrow X$ be any resolution of the singularities. Then for arbitrary ample divisors $H_1,\dots,H_{n-2}$, we have the following inequality: $$(3c_2(Y)-c_1^2(Y))\rho^*(H_1)\dots\rho^*(H_{n-2})\geq 0.$$}\par
\noindent{\bf Remark 1.12.1.} \phantomsection\label{1.12.1}An obvious corollary of the theorem above is the following: \textit{Let $X$ be a smooth three (resp. four) fold and $A$ be an ample divisor on it. Then $A\cdot c_2(X) \geq 0$ (resp. $A^2\cdot c_2(X)\geq 0$)}.

\section{Projective normality for adjoint linear series}\label{N0}
All the varieties appearing in this section are smooth. Here we will prove theorems on projective normalty and normal presentation of adjoint linear series associated to a globally generated, ample line bundle. The proofs here are  based on the philosophy that a multiplication map surjects on a variety if its restriction surjects on a certain curve. To prove this, we will use the Skoda complex defined in Section \ref{prelims} as the variety we are working on is not necessarily regular.\par
\noindent\textbf{Lemma 2.1.} \phantomsection\label{2.1}\textit{Let $X$ be a variety of dimension $n$, $n\geq 3$. Let $B$ be an ample and base point free line bundle on $X$. We further assume $h^0(B)\geq n+2$. Let $X_n$ be $X$, $X_{n-j}$ be a smooth irreducible $(n-j)$-fold chosen from the complete linear system of $|B\vert_{X_{n-j+1}}|$ (which exists by Bertini) for all $1\leq j\leq n-1$. Then the following will hold:\\
\indent (i) $H^1(K+lB\vert_{X_{n-j}})=0$ for all $0\leq j\leq n-2$, $l\geq n-1$.\\
\indent (ii) $H^0(K+nB)\otimes H^0(B)\rightarrow H^0(K+(n+1)B)$ surjects.}\par
\noindent\textit{Proof of (i).} By adjunction, $K_{X_{n-j}}=(K+jB)\vert_{X_{n-j}}$ for all $0\leq j\leq n-1$. Thus, 
by Kodaira vanishing, 
\begin{align*}
    H^1(K+lB\vert_{X_{n-j}})=H^1(K_{X_{n-j}}+(l-j)B\vert_{X_{n-j}})=0
\end{align*} 
for all $0\leq j\leq n-2$, $l\geq n-1$.\par

\noindent\textit{Proof of (ii).} Thanks to part (i), $H^0((K+lB)\vert_{X_{n-j}})\rightarrow H^0((K+lB)\vert_{X_{n-j-1}})$ surjects for all $l\geq n$, $0\leq j\leq n-2$.
We have the following diagram.\phantomsection\label{2.1.1}
\begin{equation}
\begin{tikzcd}
 0\arrow[r] & H^0(L)\otimes H^0(B\otimes \mathscr{I}) \arrow[r] \arrow[d] &  H^0(L)\otimes H^0(B) \arrow[r] \arrow[d] &  H^0(L)\otimes V \arrow[r]\arrow[d] & 0\\
0 \arrow[r] & H^0((L+B)\otimes \mathscr{I}) \arrow[r] & H^0(L+B) \arrow[r] & H^0((L+B)\vert_{X_1}) \arrow[r] & 0\tag{2.1.1}
\end{tikzcd}
\end{equation}

Here $L=K+nB$, $\mathscr{I}$ is the ideal sheaf of the curve $X_1$ in $X$ and $V$ is the cokernel of the map $H^0(B\otimes \mathscr{I})\rightarrow H^0(B)$. The bottom row is exact by part (i) and the top row is exact by the definition of $V$.\par 
Let $W$ be the vector space corresponding to the curve $X_1$ on $X$ that appears on the Skoda complex (see Definition \hyperref[1.8]{1.8}). Tensoring the following exact sequence:
\phantomsection\label{2.1.2}
\begin{equation}
    0 \longrightarrow \bigwedge\limits^{n-1}W\otimes B^{-(n-1)}\longrightarrow \dots \longrightarrow \bigwedge\limits^{2}W\otimes B^{-2} \longrightarrow W\otimes B^{-1} \longrightarrow \mathscr{I} \longrightarrow 0 \tag{2.1.2}
\end{equation}
by $L+B$, we get the following exact sequence where $L'=L+B$,\[
\begin{tikzcd}
 0 \arrow{r}{} & \bigwedge\limits^{n-1}W\otimes L'\otimes B^{-(n-1)}\arrow{r}{f_{n-1}} & \dots \arrow{r}{f_2} &  W\otimes L'\otimes B^{-1} \arrow{r}{f_1} & L'\otimes\mathscr{I} \arrow{r}{} & 0.
\end{tikzcd}
\]

To show the left most vertical map in \hyperref[2.1.1]{(2.1.1)} surjects, it is enough to prove $H^1(ker(f_1))=0$\\ as $W\subseteq H^0(B\otimes \mathscr{I})$. The following two claims prove the vanishing.\par

\noindent\textit{Claim 1\phantomsection\label{14}: $H^r(ker(f_r))=0\implies H^{r-1}(ker(f_{r-1}))=0$ for all $2\leq r\leq n-2$}.\par 
\noindent\textit{Proof:} We have the following short exact sequence:\[
\begin{tikzcd}
 0 \arrow{r}{} & ker(f_r)\arrow{r}{} & \bigwedge\limits^{r}W\otimes L'\otimes B^{-r} \arrow{r}{f_r} & ker(f_{r-1})\arrow{r}{} & 0.
\end{tikzcd}
\]
The long exact sequence of cohomology proves the claim as  $H^{r-1}(K+(n+1-r)B)=0$ since $n+1-r>0$ for $r$ in the given interval.\par

\noindent\textit{Claim 2\phantomsection\label{15}: $H^{n-2}(L'-(n-1)B)=0$.}\par 
\noindent\textit{Proof:} This is obvious from Kodaira vanishing as $H^{n-2}(L'-(n-1)B)=H^{n-2}(K+2B)=0$.\par

Thus, in order to prove the surjectivity of the middle vertical map in \hyperref[2.1.1]{(2.1.1)}, we only have to prove the surjection of the map  $H^0(L\vert_{X_1})\otimes V\rightarrow H^0((L+B)\vert_{X_1})$  as $H^0(L\vert_{X_{n-j}})\rightarrow H^0(L\vert_{X_{n-j-1}})$ already surjects for all $0\leq j\leq n-2$ by part (i).\par 
Using Lemma \hyperref[1.6]{1.6}, it is enough to prove the following inequality:
\begin{equation}
    \phantomsection\label{2.1.3}h^1((K+(n-1)B)\vert_{X_1})\leq dim(V)-2.\tag{2.1.3}
\end{equation}
\indent So, first we have to find an estimate of $dim(V)$.\par
\noindent\textit{Claim 3:\phantomsection\label{16} $h^0(B\otimes \mathscr{I})=dim(W)$}.\par 
\noindent\textit{Proof:} We tensor the exact sequence \hyperref[2.1.2]{(2.1.2)} by $B$ and get the following exact sequence:\[
\begin{tikzcd}
 0 \arrow{r}{} & \bigwedge\limits^{n-1}W\otimes B^{-(n-2)}\arrow{r}{g_{n-1}} & \dots \arrow{r}{g_3} & \bigwedge\limits^{2}W\otimes B^{-1} \arrow{r}{g_2}  & W\otimes \mathscr{O}_X \arrow{r}{g_1} & B\otimes\mathscr{I} \arrow{r}{} & 0
\end{tikzcd}
\]
\noindent So, in order to prove the claim, it is enough to show $H^0(ker(g_1))=0$ and $H^1(ker(g_1))=0$.\par 
These two vanishing can be seen from the following four facts whose proofs we omit as they are similar to \hyperref[14]{\textit{Claim 1}} and \hyperref[15]{\textit{Claim 2}}. \par 
\indent\textit{Fact 1: $H^{r-1}(ker(g_r))=0\implies H^{r-2}(ker(g_{r-1}))=0$ for all $2\leq r\leq n-2$,}\\
\indent\textit{Fact 2: $H^{n-3}(B^{-(n-2)})=0$,}\\
\indent\textit{Fact 3: $H^{r}(ker(g_r))=0\implies H^{r-1}(ker(g_{r-1}))=0$ for all $2\leq r\leq n-2$,}\\
\indent\textit{Fact 4: $H^{n-2}(B^{-(n-2)})=0$.}\par

Therefore, $dim(V)=h^0(B)-h^0(B\otimes \mathscr{I})\geq h^0(B)-(n-1)$ as $dim(W)\leq n-1$. Note that, $(K+(n-1)B)\vert_{X_1}$ is the canonical bundle of $X_1$ and consequently $h^1((K+(n-1)B)\vert_{X_1})=1$. Thus, the inequality \hyperref[2.1.3]{(2.1.3)} is verified thanks to $h^0(B)\geq n+2$.\QEDA\par


\noindent{\bf Remark 2.1.1.} Since $B$ is ample and base point free, $h^0(B)\geq n+1$. In our theorems, we are assuming that $h^0(B)\geq n+2$. Later we will give an example where $h^0(B)=4$ and $K+3B$ does not satisfy projective normality on a regular three-fold.\QEDB\par

\noindent\textbf{Lemma 2.2.} \phantomsection\label{2.2}\textit{Let $X$ be a variety of dimenson $n$, $n\geq 3$. Let $B$ be an ample and base point free line bundle on $X$. We further assume:\\
\indent (a) $K$ is nef, $K+B$ is base point free.\\
\indent (b) $h^0(K+B)\geq h^0(K)+n+1$.\\
\indent (c) $B-K$ is nef and effective divisor.\\
Let $X_n$ be $X$, $X_{n-j}$ be sufficiently general smooth irreducible $(n-j)$ fold chosen from the complete linear system of $|(K+B)\vert_{X_{n-j+1}}|$ for all $1\leq j\leq n-1$. Then the following will hold:\\
\indent (i) $H^1((2n-2)B\vert_{X_{n-j}})=0$ for all $0\leq j\leq n-2$.\\
\indent (ii) $H^0(K+(2n-1)B)\otimes H^0(K+B)\rightarrow H^0(2K+2nB)$ surjects.}\par

\noindent\textit{Proof of (i).} Adjunction gives us $K_{X_{n-j}}=((j+1)K+jB)\vert_{X_{n-j}}$ for all $0\leq j\leq n-1$. We have, 
\begin{align*}
    H^1((2n-2)B\vert_{X_{n-j}})=H^1(K_{X_{n-j}}+(((2n-2j-3)B+(j+1)(B-K))\vert_{X_{n-j}})).
\end{align*}
Note that, $2n-2j-3\geq 1$ for all $0\leq j\leq n-2$. Using Kodaira vanishing we conclude 
\begin{align*}
    H^1((2n-2)B\vert_{X_{n-j}})=0\textrm{ for all }0\leq j\leq n-2
\end{align*}
as $B-K$ is nef.\par

\noindent\textit{Proof of (ii).} Let $\mathscr{I}$ be the ideal sheaf of $X_1$ in $X$ and consequently we have $W$ as in Definition \hyperref[1.8]{1.8}. We have the following diagram \phantomsection\label{2.2.1} where $L=K+(2n-1)B$ and $V$ is the cokernel of the map $H^0((K+B)\otimes\mathscr{I})\rightarrow H^0(K+B)$:\par
\begin{equation}
    \begin{tikzcd}[row sep=large, column sep=2ex]
 0\arrow[r] & H^0(L)\otimes H^0((K+B)\otimes \mathscr{I})  \arrow[r] \arrow[d] &  H^0(L)\otimes H^0(K+B) \arrow[r] \arrow[d] &  H^0(L)\otimes V \arrow[r]\arrow[d] & 0\\
0 \arrow[r] & H^0((L+K+B)\otimes \mathscr{I}) \arrow[r] & H^0((L+K+B)) \arrow[r] & H^0((L+K+B)\vert_{X_1}) \arrow[r] & 0\tag{2.2.1}
\end{tikzcd}
\end{equation}
The bottom row is exact by Kodaira vanishing as $H^1((K+(2n-1)B)\vert_{X_{n-j}})=0$, the top row is exact by our construction. \par 
We have the following exact sequence:\[
\begin{tikzcd}
 0 \arrow{r}{} & \bigwedge\limits^{n-1}W\otimes (K+B)^{-(n-1)}\arrow{r}{} & \dots \arrow{r}{} &  W\otimes (K+B)^{-1} \arrow{r}{} & \mathscr{I} \arrow{r}{} & 0
\end{tikzcd}\]
 Tensoring by $L+K+B$ and taking cohomology, as in the proof of Lemma 2.1, we have the following two vanishings:\par 
(\textit{V-1}) $H^{r-1}(K+L+B-r(K+B))=H^{r-1}(K+(2n-2r+1)B+(r-1)(B-K))=0$ for all $2\leq r\leq n-2$ which is obvious by Kodaira vanishing since we have $B-K$ nef.\par
(\textit{V-2}) $H^{n-2}(L+K+B-(n-1)(K+B))=0$ which comes from Kodaira vanishing as well.\par

The above two vanishings show that the leftmost vertical map in \hyperref[2.2.1]{(2.2.1)} is surjective. Note that $H^0(L)\rightarrow H^0(L\vert_{X_1})$ is surjective by part (i). Consequently, by the application of Lemma \hyperref[1.6]{1.6}, we just need the following inequality:\phantomsection\label{2.2.2}
\begin{equation}
    h^1((2n-2)B\vert_{X_1})\leq dim(V)-2.\tag{2.2.2}
\end{equation}

As in the proof of \hyperref[16]{\textit{Claim 3}}, Lemma 2.1, we can see that $dim(V)\geq h^0(K+B)-(n-1)$. Still, we have to estimate $h^1((2n-2)B\vert_{X_1})$.
We have the short exact sequence:
\[
\begin{tikzcd}
 0 \arrow{r}{} & (-K-B)\vert_{X_2} \arrow{r}{} \arrow{r}{} & \mathscr{O}_{X_2} \arrow{r}{} & \mathscr{O}_{X_1} \arrow{r}{} & 0.
\end{tikzcd}
\]
Tensoring this by $(2n-2)B$ gives: 
\[
\begin{tikzcd}
 0 \arrow{r}{} & (-K+(2n-3)B)\vert_{X_2} \arrow{r}{} \arrow{r}{} & (2n-2)B\vert_{X_2} \arrow{r}{} & (2n-2)B\vert_{X_1} \arrow{r}{} & 0.
\end{tikzcd}
\]
Consequently, we have the long exact sequence:
\[
\begin{tikzcd}
  \dots \arrow{r}{} & H^1((2n-2)B\vert_{X_2}) \arrow{r}{} \arrow{r}{} & H^1((2n-2)B\vert_{X_1}) \arrow{r}{} & H^2((-K+(2n-3)B)\vert_{X_2}) \arrow{r}{} & \dots
\end{tikzcd}
\]
Since $H^1((2n-2)B\vert_{X_2})=0$ by (i), we get $h^1((2n-2)B\vert_{X_1})\leq h^2((-K+(2n-3)B)\vert_{X_2})$.
Now, we have, $h^2((-K+(2n-3)B)\vert_{X_2})$\hspace{3pt}  $\myeq$ \hspace{3pt} $h^0(((n-1)K+(n-2)B+K-(2n-3)B)\vert_{X_2})=h^0(K\vert_{X_2}-(n-1)(B-K)\vert_{X_2})$.\par 
Note that, assumption (c) gives us $h^0(K\vert_{X_2}-(n-1)(B-K))\vert_{X_2})\leq h^0(K\vert_{X_2})$.
The long exact sequence associated to the following short exact sequence: 
\[
\begin{tikzcd}
 0 \arrow{r}{} & (-B)\vert_{X_{n-j+1}} \arrow{r}{} \arrow{r}{} & K\vert_{X_{n-j+1}} \arrow{r}{} & K_{X_{n-j}} \arrow{r}{} & 0
\end{tikzcd}
\]
shows us (by Kodaira vanishing) that $h^0(K\vert_{X_{n-j}})=h^0(K\vert_{X_{n-j+1}})$ for all $0\leq j\leq n-2$. Consequently, $h^0(K\vert_{X_2})=h^0(K)$.
Thus, to show inequality \hyperref[2.2.2]{(2.2.2)} it is enough to show $h^0(K)\leq h^0(K+B)-(n+1)$ which we have, thanks to assumption (b).\QEDA\par

\noindent{\bf Remark 2.2.1.} \phantomsection\label{2.2.01}We always have $h^0(K+B)\geq h^0(K)+n$ on any $n$-fold if $K+B$ and $B$ are ample and base point free.\par 
\noindent\textit{Proof.} Note that, $h^0(K+B)-h^0(K)$ is the dimension of the cokernel $V$ of the map $H^0(K)\rightarrow H^0(K+B)$ in $H^0((K+B)\vert_{X_{n-1}})$ where $X_{n-1}$ is a smooth irreducible divisor chosen from the complete linear system of $B$. Notice that $V$ is the linear subsystem of the complete linear series $|(K+B)\vert_{X_{n-1}}|$ obtained by pulling back the base point free complete linear series $|K+B|$ on $X$ by the embedding $X_{n-1}\hookrightarrow X$. Consequently $|V|$ is globally generated and $dim(V)\geq n$ as the morphism induced by $|V|$ is the composite of the embedding $i$ and a finite morphism (given on $X$ by $|K+B|$) and is hence finite on an $n-1$ dimensional variety $X_{n-1}$.\QEDB\par

\noindent{\bf Remark 2.2.2.} \phantomsection\label{2.2.02}Let $X$ be a variety of dimension $n$ with nef canonical bundle $K$. Let $B$ be an ample and base point free line bundle such that $B+K$ is globally generated, $h^0(B)\geq n+2$ and $H^1(B)=0$. Then $h^0(K+B)\geq h^0(K)+n+1$. \par 
\noindent\textit{Proof.} The assertion is trivial if $h^0(K)=0$ or if $K=\mathscr{O}_X$.
Otherwise, we have the short exact sequence:
\[
\begin{tikzcd}
 0 \arrow{r}{} & B \arrow{r}{} & B+K \arrow{r}{} & (B+K)\vert_{\mathscr{K}} \arrow{r}{} & 0
\end{tikzcd}
\]
where $\mathscr{K}$ is a non zero effective divisor chosen from the linear system of $K$.\par 
From the long exact sequence, we get that $h^0(K+B)=h^0(B)+h^0((B+K)\vert_{\mathscr{K}})$.
But $h^0((B+K)\vert_{\mathscr{K}})\geq h^0(K\vert_{\mathscr{K}})$. Thus,
$h^0(K+B)=h^0(B)+h^0((B+K)\vert_{\mathscr{K}})\geq n+2+h^0(K)-1$.\QEDB\par

\noindent{\bf Remark 2.2.3.} \phantomsection\label{2.2.03}Let $X$ be a variety of dimension $n$ with nef canonical bundle $K$ and $H^{n-1}(\mathscr{O}_X)=0$. Let $B$ be an ample and base point free line bundle on $X$ such that $B+K$ is globally generated and $h^0(B)\geq n+2$. Then $h^0(K+B)\geq h^0(K)+n+1$.\par 
\noindent\textit{Proof.} Again, we can assume that $K$ is a non zero effective divisor.
The long exact sequence associated to the short exact sequence: 
\[
\begin{tikzcd}
 0 \arrow{r}{} & K \arrow{r}{} & B+K \arrow{r}{} & (B+K)\vert_{\mathscr{B}} \arrow{r}{} & 0
\end{tikzcd}
\]
gives $h^0(K+B)=h^0(K)+h^0((B+K)\vert_{\mathscr{B}})$
(here $\mathscr{B}$ is a sufficiently general non zero effective divisor chosen from the linear system of $B$).
Now we have $h^0((B+K)\vert_{\mathscr{B}})\geq h^0(B\vert_{\mathscr{B}})$.
Hence $h^0(K+B)\geq h^0(K)+h^0(B\vert_{\mathscr{B}})\geq h^0(K)+n+1$.\QEDB\par

Now we prove our first main result that gives the projective normality of $K+nB$ on a regular $n$ dimensional variety under some assumptions.\par

\noindent{\textbf{Theorem 2.3.}} \phantomsection\label{2.3}\textit{Let $X$ be a variety of dimension $n$, $n\geq 3$. Let $B$ be an ample and base point free line bundle on $X$. We further assume:\\
\indent (a) $K$ is nef, $K+B$ is base point free.\\
\indent (b) $h^0(B)\geq n+2$.\\
\indent (c) $h^0(K+B)\geq h^0(K)+n+1$.\\
\indent (d) $B-K$ is nef and effective.\\
Then $K+lB$ is very ample and it embeds $X$ as a projectively normal variety for all $l\geq n$.}\par

\noindent\textit{Proof.} We need to prove $H^0((K+lB)^{\otimes k})\otimes H^0(K+lB)\rightarrow H^0((K+lB)^{\otimes k+1})$ surjects $\forall k\geq 1$.\par

\textit{Step 1: \phantomsection\label{1}$H^0(k(K+lB)+rB)\otimes H^0(B)\rightarrow H^0(k(K+lB)+(r+1)B)$ surjects for $k\geq 2$, $l\geq n$, $r\geq 0$.}\\
This comes from CM lemma (Lemma \hyperref[1.7]{1.7}) once we note that,
\begin{equation*}
    H^i(k(K+lB)+(r-i)B)=H^i(K+(k-1)K+rB+(kl-i)B)=0\textrm{ for all }1\leq i\leq n
\end{equation*}
by Kodaira vanishing.\par

\textit{Step 2: \phantomsection\label{2}$H^0(k(K+lB)+(l-1)B)\otimes H^0(K+B)\rightarrow H^0((k+1)K+(kl+l)B)$ surjects for $k\geq 2$, $l\geq n$.}
This again comes from CM lemma (Lemma \hyperref[1.7]{1.7}). Indeed,
\begin{equation*}
    H^i(k(K+lB)+(l-1)B-iK-iB)=H^i(K+kK+(kl+l-1-i)B-(1+i)K).
\end{equation*}
But $kl+l-1-i\geq 3n-1-i\geq  2n-1> 1+n$ for all $1\leq i\leq n$. 
Since assumption (d) shows us that $B-K$ is nef, Kodaira vanishing gives 
\begin{equation*}
    H^i(k(K+lB)+(l-1)B-iK-iB)=0\textrm{ for all }1\leq i\leq n.
\end{equation*}\par

\textit{Step 3: $H^0((K+lB)^{\otimes k})\otimes H^0(K+lB)\rightarrow H^0((k+lB)^{\otimes k+1})$ surjects for $k\geq 2$, $l\geq n$, $r\geq 0$.}
This comes from \textit{Steps \hyperref[1]{1}, \hyperref[2]{2}} and the Observation \hyperref[1.2]{1.2}.\par

So, we only need to prove $H^0(K+lB)\otimes H^0(K+lB)\rightarrow H^0((K+lB)^{\otimes 2})$ surjects for all $l\geq n$.\par

\textit{Step 4: \phantomsection\label{4}$H^0(K+lB)\otimes H^0(B)\rightarrow H^0(K+(l+1)B)$ surjects for $l>n$.}
This comes from CM lemma (Lemma \hyperref[1.7]{1.7}) once we note that $H^i(K+(l-i)B)=0$ for all $1\leq i\leq n$ by Kodaira vanishing since $l-i>0$.\par

\textit{Step 5: \phantomsection\label{5}$H^0(K+(2l-1)B)\otimes H^0(K+B)\rightarrow H^0(2K+2lB)$ surjects for $l>n$.} To see this, first
note that $H^i(K+(2l-1)B-iK-iB)=H^i(K+(2l-i-1)B-iK)$.
Now, $2l-i-1\geq2n-i+1> i$, thanks to $l\geq n+1$ and $1\leq i\leq n$.
Since $B-K$ is nef, Kodaira vanishing implies $H^i(K+(2l-1)B-iK-iB)=0$.
Hence by CM lemma (Lemma \hyperref[1.7]{1.7}), we are done.\par

\textit{Step 6: \phantomsection\label{6}$H^0(K+lB)\otimes H^0(K+lB)\rightarrow H^0(2K+2lB)$ surjects for $l>n$.}
This comes from \textit{Steps \hyperref[4]{4}, \hyperref[5]{5}} and the Observation \hyperref[1.2]{1.2}.\par

So, we only need to prove $H^0(K+nB)\otimes H^0(K+nB)\rightarrow H^0(2K+2nB)$ surjects which is our final step. \par

\textit{Step 7: \phantomsection\label{7}$H^0(K+nB)\otimes H^0(K+nB)\rightarrow H^0(2K+2nB)$ surjects.} This is because
in Lemma \hyperref[2.1]{2.1}, we have already proved $H^0(K+nB)\otimes H^0(B)\rightarrow H^0(K+(n+1)B)$ surjects and in \textit{Step \hyperref[4]{4}} we have showed $H^0(K+lB)\otimes H^0(B)\rightarrow H^0(K+(l+1)B)$ surjects for $l>n$.
Using Observation \hyperref[1.2]{1.2}, we only need to show the surjection of $H^0(K+(2n-1)B)\otimes H^0(K+B)\rightarrow H^0(2K+2nB)$ which we have proved in Lemma \hyperref[2.2]{2.2}.\QEDA\par

\noindent{\bf Remark 2.3.1.} \phantomsection\label{2.3.01}Let $X$ be a $n$ dimensional variety. Let $B$ be a globally generated, ample line bundle on $X$. We further assume that $B-K$ is a non-zero effective divisor. If $p_g(X)\geq 2$, then $h^0(B)\geq n+p_g$.\par 
\noindent \textit{Proof.} The long exact sequence associated to the short exact sequence: 
\[
\begin{tikzcd}
 0 \arrow{r}{} & K \arrow{r}{} & B \arrow{r}{} & B\vert_D \arrow{r}{} & 0
\end{tikzcd}
\]
shows that the cokernel of the map $H^0(K)\rightarrow H^0(B)$ is a base point free linear subsystem of the base point free complete linear system of $H^0(B\vert_D)$ ($D$ is an element of the linear series of $B-K$). By the same argument used in the proof of Remark \hyperref[2.2.01]{2.2.1}, we have $h^0(B)-p_g\geq n$.\QEDB\par

Remark \hyperref[2.3.01]{2.3.1} and Remark \hyperref[2.2.02]{2.2.2}, \hyperref[2.2.03]{2.2.3} allow us to deduce a corollary of Theorem \hyperref[2.3]{2.3} which we state below.\par

\noindent\textit{{\bf Corollary 2.4.} \phantomsection\label{2.4}Let $X$ be a variety of dimension $n\geq 3$ with $p_g\geq 2$. Let $B$ be an ample, globally generated line bundle on $X$. Assume $K$ is nef and $B-K$ is a nef, non-zero, effective divisor. Further assume that $B+K$ is globally generated. If either $H^1(B)=0$ or $H^{n-1}(\mathscr{O}_X)=0$ then $K+nB$ will be very ample and it will embed $X$ as a projectively normal variety.}\QEDA\par

Now we produce examples to discuss the sharpness of our conditions. In our first example, we construct a regular variety of general type and an ample, globally generated line bundle $B$ on it that satisfies all the conditions of Theorem \hyperref[2.3]{2.3} except the condition (b) and show that the line bundle $K+nB$ does not satisfy the property $N_0$.\par

\noindent\textbf{Example 2.5.} \phantomsection\label{2.5}Consider a double cover $X$ of $\mathbb{P}^{n+1}$ ramified along an $n$-fold of degree $2n+4$, $n\geq 3$. Let the natural finite morphism from $X$ to $\mathbb{P}^{n+1}$ be $f$. The unique line bundle associated to this cover is $\mathscr{O}(n+2)$. We have $f_*(\mathscr{O}_X)= \mathscr{O}\oplus \mathscr{O}(-n-2)$, $H^1(\mathscr{O}_X)=0$ and $K_X=\mathscr{O}_X$.\par 
Consider $B=f^*\mathscr{O}(1)$. Clearly $B$ is ample and base point free,
\begin{equation*}
    H^0(B)=H^0(f_*(B))=H^0(\mathscr{O}(1)\oplus \mathscr{O}(-n-1)) \implies h^0(B)=n+2.
\end{equation*}
 Kodaira vanishing shows that $H^1(rB)=0$ for all $r\geq 1$. \par
  Let $Y \in |B|$ be smooth irreducible $n$-fold given by Bertini's Theorem. Consider the line bundle $B|_Y$ on $Y$. This is again ample and base point free and by adjunction $K_Y=B|_Y$. So $Y$ is a smooth $n$-fold of general type. Consider the following exact sequence:\phantomsection\label{2.5.1}
\begin{equation}
    \begin{tikzcd}
0 \arrow[r] & B^* \arrow[r] & \mathscr{O}_X \arrow[r] & \mathscr{O}_Y \arrow[r] & 0. \tag{2.5.1}
\end{tikzcd}
\end{equation}  

Taking cohomology and using Kodaira vanishing, we get  $H^1(\mathscr{O}_Y)=0$. Tensoring the sequence \hyperref[2.5.1]{(2.5.1)} by $B$ and taking cohomology gives us $h^0(B|_Y)=h^0(B)-1=n+1$. Hence $B|_Y$ does not satisfy the condition (b) of Theorem \hyperref[2.3]{2.3}.\par 
We have that $K_Y+B|_Y=2B|_Y$ and is hence base point free. Clearly $K_Y=B|_Y$ is nef since it is ample. Also $B|_Y-K_Y=\mathscr{O}_Y$ and is hence nef and effective. Now we show that,
\begin{equation*}
    h^0(K_Y+B|_Y)\geq h^0(K_Y)+n+1\textrm{ i.e. }h^0(2B|_Y)\geq h^0(B|_Y)+n+1.
\end{equation*}
We have that $H^0(2B)=H^0(f^*\mathscr{O}(2))=H^0(\mathscr{O}(2)\oplus \mathscr{O}(-n)) \implies h^0(2B)=h^0(\mathscr{O}(2))=\binom{n+2}{2}+(n+2)$.\par 
Tensoring the exact sequence \hyperref[2.5.1]{(2.5.1)} by $2B$ and taking the cohomology shows that 
\begin{equation*}
    h^0(2B|_Y)= h^0(2B)-h^0(B)=\binom{n+2}{2}.
\end{equation*}
Now, $h^0(B|_Y)+n+1=2n+2$. Since $n\geq 3$ we have that $h^0(K_Y+B|_Y)\geq h^0(K_Y)+n+1$. 
We have showed that $B|_Y$ satisfies all the conditions in Theorem \hyperref[2.3]{2.3} except (b). Now we prove that $K_Y+nB|_Y=(n+1)B|_Y$ does not satisfy property $N_0$.\par

We have that $K_Y+nB|_Y=(n+1)B\vert_Y$. Suppose $(n+1)B|_Y$ satisfies the property $N_0$. 
Hence for a curve section $C \in |B_Y|$ we have that $(n+1)B|_C$ is very ample. We also have that $K_C=nB|_C$ and hence $(n+1)B|_C=K_C+B|_C$. \par 
Now $deg(B|_C)=B^{n+1}=2H^{n+1}=2$ where $H$ is a hyperplane section of $\mathbb{P}^{n+1}$ since the map $f$ is 2-1. But $K_C+E$ cannot be very ample if $E$ is an effective divisor of degree $2$.\QEDB\par

Now we give an example of a variety and an ample, globally generated line bundle $B$ for which $K+nB$ does not satisfy the property $N_0$, where $B-K$ is neither nef nor effective although the geometric genus of the variety is large (see Corollary \hyperref[2.4]{2.4}).\par

\noindent\textbf{Example 2.6.} \phantomsection\label{2.6}Consider $X$ a cyclic double cover of $\mathbb{P}^n$ ramified along hypersurface of degree $2r$. Denote by $f$ the natural morphism from $X$ to $\mathbb{P}^n $. Let $B=f^*(\mathscr{O}(1))$. We have that,\par  \noindent $f_*(\mathscr{O}_X)=\mathscr{O}\oplus \mathscr{O}(-r)$, $K_X=f^*(\mathscr{O}(-n-1+r))$, $K_X+B=f^*(\mathscr{O}(-n+r))$, $B-K_X=f^*(\mathscr{O}(n+2-r))$.\par 
We can see that for $r\geq n+3$, $B-K_X$ is not nef. However by making $r$ large enough we can make $p_g$ as large as we wish to and in particular make $p_g\geq 2$. We also have $H^1(B)=0$. We now show that for $r\geq n+3 $, $K_X+nB$ is not projectively normal. Indeed,
\begin{equation*}
    K_X+nB=f^*(\mathscr{O}(r-1))\implies H^0(K+nB)=H^0(\mathscr{O}(r-1)\oplus \mathscr{O}(-1))=H^0(\mathscr{O}(r-1)).
\end{equation*}
Now $H^0(2K_X+2nB)=H^0(f^*(\mathscr{O}(2r-2)))=H^0(\mathscr{O}(2r-2))\oplus H^0(\mathscr{O}(r-2))$. If $r\geq 2$ we can clearly see that $K+nB$ is not projectively normal. Hence we can see that the condition $B-K_X$ nef and effective is essential in Corollary \hyperref[2.4]{2.4}.\QEDB

\section{Normal presentation for adjoint linear series}\label{N1}
Our goal is to prove results concerning the $N_1$ property of adjoint bundles. Unlike the previous section, first we prove results for regular varieties and then we prove a weaker result for irregular varieties.  We prove three technical lemmas to begin with. The proofs are again based on the same philosophy that a multiplication map surjects on a variety iff it surjects on a curve section. All the varieties in this section are smooth.\par

\noindent\textbf{Lemma 3.1.} \phantomsection\label{3.1}\textit{Let $X$ be a regular variety of dimension $n$, $n\geq 3$. Let $B$ be an ample and base point free line bundle on $X$. We further assume $h^0(B)\geq n+2$. 
Let $X_n$ be $X$ and let $X_{n-j}$ be a smooth irreducible $(n-j)$-fold chosen from the complete linear system $|B\vert_{X_{n-j+1}}|$ (which exists by Bertini) for all $1\leq j\leq n-1$.
Then the map 
\begin{equation*}
    H^0((K+nB)\vert_{X_{n-j}})\otimes H^0(B\vert_{X_{n-j}})\rightarrow H^0((K+(n+1)B)\vert_{X_{n-j}})
\end{equation*}
surjects for $0\leq j\leq n-1$.}\par

\noindent\textit{Proof.} To start with, notice that $X_{n-j}$ is regular for all $j$. Indeed, it can easily be seen by taking cohomology of the exact sequence
\begin{equation*}
\begin{tikzcd}
0\arrow[r] & -B\vert_{X_{n-j+1}}\arrow[r] &\mathscr{O}_{X_{n-j+1}}\arrow[r] &\mathscr{O}_{X_{n-j}}\arrow[r] & 0.
\end{tikzcd}
\end{equation*}
Because of the vanishing Lemma \hyperref[2.1]{2.1 (i)}, by the repeated application of Lemma \hyperref[1.3]{1.3}, it is enough to prove $H^0((K+nB)\vert_{X_1})\otimes H^0(B\vert_{X_1})\rightarrow H^0((K+(n+1)B)\vert_{X_1})$ surjects.
To show this using Lemma \hyperref[1.6]{1.6}, we have to prove the inequality $h^1((K+(n-1)B)\vert_{X_1})\leq h^0(B\vert_{X_1})-2$ which follows directly from our assumption that $h^0(B)\geq n+2$.\QEDA\par

\noindent\textit{{\bf Lemma 3.2.} \phantomsection\label{3.2}Let $X$ be a regular $n$-fold, $n\geq 3$. Let $B$ be an ample and base point free line bundle on $X$. We further assume:\\
\indent (a) $K$ is nef, $K+B$ is base point free.\\
\indent (b) $h^0(K+B)\geq h^0(K)+n+1$.\\
\indent (c) $(n-2)B-(n-1)K$ is nef and effective.\\
Let $X_n$ be $X$ and let $X_{n-j}$ be a sufficiently general smooth irreducible $(n-j)$-fold chosen from the complete linear system $|(K+B)\vert_{X_{n-j+1}}|$ for all $1\leq j\leq n-1$.
Then the following will hold:\\
\indent(i) $H^1((2n-3)B\vert_{X_{n-j}})=0$ for all $0\leq j\leq n-2$.\\
\indent (ii) $H^0((K+(2n-2)B)\vert_{X_{n-j}})\otimes H^0((K+B)\vert_{X_{n-j}})\rightarrow H^0((2K+(2n-1)B)\vert_{X_{n-j}})$ surjects for all $0\leq j\leq n-1$.}\par

\noindent\textit{Proof of (i).} Adjunction gives us $K_{X_{n-j}}=((j+1)K+jB)\vert_{X_{n-j}}$ for all $0\leq j\leq n-1$.
We have 
\begin{equation*}
    H^1((2n-3)B\vert_{X_{n-j}})=H^1(K_{X_{n-j}}+(B+(2n-4-j)B-(j+1)K)\vert_{X_{n-j}}).
\end{equation*}
Note that, 
\begin{equation*}
    B-\displaystyle\frac{j+1}{2n-4-j}K=B-\displaystyle\frac{n-1}{n-2}K+\displaystyle\frac{n-1}{n-2}K-\displaystyle\frac{j+1}{2n-4-j}K.
\end{equation*}
We have $n-1\geq j+1$ and $n-2\leq 2n-4-j$ for all $0\leq j\leq n-2$. Consequently, $(2n-4-j)B-(j+1)K$ is nef as $K$ and $B-\displaystyle\frac{n-1}{n-2}K$ are nef.
Using Kodaira vanishing we conclude $H^1((2n-3)B\vert_{X_{n-j}})=0$ for all $0\leq j\leq n-2$.\par

\noindent\textit{Proof of (ii).} As in the proof of Lemma \hyperref[3.1]{3.1}, $X_{n-j}$ is regular forr all $j$. Repeated application of Lemma \hyperref[1.3]{1.3} shows that it is enough to prove the lemma for $j=n-1$. Hence, we have to prove the surjection of 
\begin{equation*}
    H^0((K+(2n-2)B)\vert_{X_1})\otimes H^0((K+B)\vert_{X_1})\rightarrow H^0((2K+(2n-1)B)\vert_{X_1}).
\end{equation*}
Application of Lemma \hyperref[1.6]{1.6} shows us it is enough to check the following inequality: \phantomsection\label{3.2.1}
\begin{equation}
    h^1((2n-3)B\vert_{X_1})\leq h^0((K+B)\vert_{X_1})-2.\tag{3.2.1}
\end{equation}
We have the short exact sequence:
\[
\begin{tikzcd}
 0 \arrow{r}{} & (-K-B)\vert_{X_2} \arrow{r}{} \arrow{r}{} & \mathscr{O}_{X_2} \arrow{r}{} & \mathscr{O}_{X_1} \arrow{r}{} & 0.
\end{tikzcd}
\]
Tensoring this by $(2n-3)B$ gives: 
\[
\begin{tikzcd}
 0 \arrow{r}{} & (-K+(2n-4)B)\vert_{X_2} \arrow{r}{} \arrow{r}{} & (2n-3)B\vert_{X_2} \arrow{r}{} & (2n-3)B\vert_{X_1} \arrow{r}{} & 0.
\end{tikzcd}
\]
Consequently, we have the long exact sequence:
\[
\begin{tikzcd}
  \dots \arrow{r}{} & H^1((2n-3)B\vert_{X_2}) \arrow{r}{} \arrow{r}{} & H^1((2n-3)B\vert_{X_1}) \arrow{r}{} & H^2((-K+(2n-4)B)\vert_{X_2}) \arrow{r}{} & \dots
\end{tikzcd}
\]
Since $H^1((2n-3)B\vert_{X_2})=0$ by (i), we get $h^1((2n-3)B\vert_{X_1})\leq h^2((-K+(2n-4)B)\vert_{X_2})$. We have,
\begin{align*}
    h^2((-K+(2n-4)B)\vert_{X_2}) & \myeq h^0(((n-1)K+(n-2)B+K-(2n-4)B)\vert_{X_2})\\
     & \;\; = h^0((K+(n-1)K-(n-2)B)\vert_{X_2}).
\end{align*}
\indent Note that, assumption (c) gives us $h^0((K+(n-1)K-(n-2)B)\vert_{X_2})\leq h^0(K\vert_{X_2})$. 
The long exact sequence associated to the following short exact sequence 
\[
\begin{tikzcd}
 0 \arrow{r}{} & (-B)\vert_{X_{n-j+1}} \arrow{r}{} \arrow{r}{} & K\vert_{X_{n-j+1}} \arrow{r}{} & K\vert_{X_{n-j}} \arrow{r}{} & 0
\end{tikzcd}
\]
shows us (by Kodaira vanishing) that $h^0(K\vert_{X_{n-j}})=h^0(K\vert_{X_{n-j+1}})$ for all $0\leq j\leq n-2$. Consequently we get that $h^0(K\vert_{X_2})=h^0(K)$.\par 
Thus, in order to show \hyperref[3.2.1]{(3.2.1)} it is enough to show $h^0(K)\leq h^0((K+B)\vert_{X_1})-2$ which comes from assumption (b). Indeed, tensoring the above exact sequence by $B\vert_{X_{n-j+1}}$ and taking cohomology (recall that $X_{n-j+1}$ is regular), one sees easily that $h^0((K+B)\vert_{X_{n-j}})=h^0((K+B)\vert_{X_{n-j+1}})-1$.\QEDA\par

\noindent\textit{{\bf{Lemma 3.3.}} \phantomsection\label{3.3}Let $X$ be a regular $n$-fold, $n\geq 3$. Let $B$ be an ample and base point free line bundle on $X$. We further assume:\\
\indent (a) $K$ is nef, $K+B$ is base point free.\\
\indent (b) $h^0(B)\geq n+2$.\\
\indent (c) $h^0(K+B)\geq h^0(K)+n+1$.\\
\indent (d) $(n-2)B-(n-1)K$ is nef and non-zero effective divisor.\\
Let $X_n$ be $X$ and let $X_{n-j}$ be a sufficiently general smooth irreducible $(n-j)$ fold chosen from the complete linear system of $|B\vert_{X_{n-j+1}}|$ for all $1\leq j\leq n-1$. Let $L$ be $K+lB$ where $l\geq n$. Then} 
\begin{equation*}
\begin{tikzcd}
H^0(M_{L\vert_{X_{n-j}}}\otimes L\vert_{X_{n-j}})\otimes H^0(B\vert_{X_{n-j}})\longrightarrow H^0(M_{L\vert_{X_{n-j}}}\otimes L\vert_{X_{n-j}}\otimes B\vert_{X_{n-j}} )
\end{tikzcd}
\end{equation*}
\textit{surjects for all $0\leq j\leq n-1$.}\par

\noindent\textit{Proof.} For simplicity, we only prove the assertion for $l=n$ that is for $L=K+nB$. The proof for $l>n$ is similar. We prove the case for $L=K+nB$ by induction on $j$. Before starting the induction, we first prove the following claim.\par

\noindent \textit{Claim:} In the context of our theorem we have $B^n \geq 4$.\par
\noindent \textit{Proof of the Claim:} Let $h^0(B)=r+1$. Let $f$ be the finite morphism induced by the ample and base point free line bundle $B$.
We have that $B^n=deg(f)\cdot deg(Y)$ where $Y$ is the scheme theoretic image. Now, the codimension of $Y$ in $\mathbb{P}^r \geq 1$ and hence $deg(Y) \geq 2$. So the only way $B^n < 4$ is when the following happens.\par
 \textit{Case 1:} $deg(f)=1$ and $deg(Y)=2$ and hence $codim(Y)=1$.
 In this case we have that $Y$ is a variety of minimal degree and it is either a smooth quadric hypersurface or a cone over a smooth rational normal scroll or a cone over the Veronese embedding of $\mathbb{P}^2$ (see \cite{Ei}). In all three cases $Y$ is normal. Indeed, the first case is trivial. The second and third case follows from the fact that a cone over a projectively normal embedding is normal. Now $f$ is a finite birational map between normal varieties and is hence an isomorphism. Consequently, the image is a smooth rational normal scroll whose canonical divisor is negative ample. This contradicts the hypothesis (a).\par
  \textit{Case 2:} $deg(f)=1$ and $deg(Y)=3$ and $codim(Y)=2$.
 In this case again $Y$ is a variety of minimal degree and hence a normal variety and we have that $f$ is an isomorphism which leads to a contradiction as before.\par
  \textit{Case 3:} $deg(f)=1$ and $deg(Y)=3$ and $codim(Y)=1$.
 In this case consider a general curve section $C$ of $|B|$ in $X$. It is the pullback of a general curve section $D$ of $\mathscr{O}(1)$ in $Y$. By Bertini we have that $C$ can be taken to be smooth and irreducible and since $f$ is surjective we have that $D$ is reduced and irreducible. Notice that $D$ is a plane curve since the codimension of $Y$ was 1. Moreover, $D$ is a plane curve of degree $3$ and hence we have that $p_a(D)=1$. $C$ is the normalization of $D$ and hence $g(C) \leq 1$. But we have that $2g(C)-2=(n-1)B^n+B^{n-1}K_C$ and hence $g(C) \geq 4$ since $B^n=3$, $n \geq 3$ and $K_C$ is nef. So we have a contradiction.\par

Now we start our induction on $j$. We aim to show the following. \par

\noindent\textit{Induction Step: Suppose $H^0(M_{L\vert_{X_{n-j}}}\otimes L\vert_{X_{n-j}})\otimes H^0(B\vert_{X_{n-j}})\rightarrow H^0(M_{L\vert_{X_{n-j}}}\otimes L\vert_{X_{n-j}}\otimes B\vert_{X_{n-j}} )$ surjects for some $1\leq j\leq n-1$. Then 
\begin{equation*}
    H^0(M_{L\vert_{X_{n-j+1}}}\otimes L\vert_{X_{n-j+1}})\otimes H^0(B\vert_{X_{n-j+1}})\rightarrow H^0(M_{L\vert_{X_{n-j+1}}}\otimes L\vert_{X_{n-j+1}}\otimes B\vert_{X_{n-j+1}} )
\end{equation*}
surjects.}\par
\noindent\textit{Proof of Induction Step:} First we prove \phantomsection\label{3.3.1}
\begin{equation}
    H^1(M_{L\vert_{X_{n-j+1}}}\otimes L\vert_{X_{n-j+1}}\otimes (B\vert_{X_{n-j+1}})^*)=0\textrm{ for $ 1\leq j\leq n-1$}.\tag{3.3.1}
\end{equation}
We have the short exact sequence:
\[
\begin{tikzcd}
 0 \arrow{r}{} & M_{L\vert_{X_{n-j+1}}}\otimes L'\vert_{X_{n-j+1}} \arrow{r}{}  & H^0(L\vert_{X_{n-j+1}})\otimes L'\vert_{X_{n-j+1}}  \arrow{r}{} & (L+L')\vert_{X_{n-j+1}} \arrow{r}{} & 0
\end{tikzcd}
\]
where $L'=K+(n-1)B$. In order to prove  \hyperref[3.3.1]{(3.3.1)}, it is enough to prove
\begin{equation*}
    H^0(L\vert_{X_{n-j+1}})\otimes H^0(L'\vert_{X_{n-j+1}}) \rightarrow H^0((L+L')\vert_{X_{n-j+1}})
\end{equation*}
surjects since according to Lemma \hyperref[2.1]{2.1 (i)}, $H^1(L'\vert_{X_{n-j+1}})=0$.\par 

To prove this surjection with the help of Observation \hyperref[1.2]{1.2}, we need to prove the following:
\begin{equation}
     H^0((K+lB)\vert_{X_{n-j+1}})\otimes H^0(B\vert_{X_{n-j+1}}) \rightarrow H^0((L+(l+1)B)\vert_{X_{n-j+1}}) \textrm{ surjects for all }l>n.\phantomsection\label{3.3.2}\tag{3.3.2}
\end{equation}
\begin{equation}
    H^0((K+nB)\vert_{X_{n-j+1}})\otimes H^0(B\vert_{X_{n-j+1}}) \rightarrow H^0((L+(n+1)B)\vert_{X_{n-j+1}})\textrm{ surjects. }\phantomsection\label{3.3.3} \tag{3.3.3}
\end{equation}
\begin{equation}
    H^0((K+(2n-2)B)\vert_{X_{n-j+1}})\otimes H^0((K+B)\vert_{X_{n-j+1}}) \rightarrow H^0((2K+(2n-1)B)\vert_{X_{n-j+1}})\textrm{ surjects.}\phantomsection\label{3.3.4}\tag{3.3.4}
\end{equation}

We use CM Lemma (Lemma \hyperref[1.7]{1.7}) to prove \hyperref[3.3.2]{(3.3.2)}. Recall that $K_{X_{n-j+1}}=(K+(j-1)B)\vert_{X_{n-j+1}}$.
Hence by Kodaira vanishing,
\begin{equation*}
    H^i((K+(l-i)B)\vert_{X_{n-j+1}})=H^i((K+(j-1)B)\vert_{X_{n-j+1}}+((l-i-j+1)B)\vert_{X_{n-j+1}})=0
\end{equation*}
for all $1\leq i\leq n-j+1$.  \par

We have already proved \hyperref[3.3.3]{(3.3.3)} in Lemma \hyperref[3.1]{3.1}.\par

For simplicity we do some re-indexing to prove \hyperref[3.3.4]{(3.3.4)} only. We will show that 
\begin{equation*}
    H^0((K+(2n-2)B)\vert_{X_{n-j}})\otimes H^0((K+B)\vert_{X_{n-j}}) \rightarrow H^0((2K+(2n-1)B)\vert_{X_{n-j}})
\end{equation*}
surjects for all $0\leq j\leq n-2$. We have already proved the surjection when $j=0$ in Lemma \hyperref[3.2]{3.2}. So, we assume $1\leq j\leq n-2$.
Our obvious choice is to use the CM Lemma (Lemma \hyperref[1.7]{1.7}).\par 
For all $1\leq i\leq n-j-1$, the following holds:
\begin{equation*}
    H^i(((1-i)K+(2n-2-i)B)\vert_{X_{n-j}})=H^i((K+jB+(2n-2-2i-j)B+i(B-K))\vert_{X_{n-j}})=0
\end{equation*}
as $B-K$ is nef and $2n-2-2i-j>0$ for $i$ in the given range.\par 
Notice that, $H^{n-j}(((n+j-2)B-(n-j-1)K)\vert_{X_{n-j}})$  $\myeq$  $H^0(((n-1)K-(n-2)B-(j-1)K)\vert_{X_{n-j}})$.
Now, $((n-1)K-(n-2)B-(j-1)K)\vert_{X_{n-j}}$ is negative nef and $(n-1)K-(n-2)B$ is negative of a non-zero effective divisor and consequently $H^{n-j}(((n+j-2)B-(n-j-1)K)\vert_{X_{n-j}})=0$.\par

Since we have proved \hyperref[3.3.1]{(3.3.1)}, to finish the proof of \textit{Induction Step} using Lemma \hyperref[1.3]{1.3}, it is enough to prove the following map
\begin{equation*}
    H^0(M_{L\vert_{X_{n-j+1}}}\otimes L\vert_{X_{n-j+1}}\otimes \mathscr{O}_{X_{n-j}})\otimes H^0(B\vert_{X_{n-j+1}}\otimes\mathscr{O}_{X_{n-j}})\rightarrow H^0(M_{L\vert_{X_{n-j+1}}}\otimes L\vert_{X_{n-j+1}}\otimes B\vert_{X_{n-j+1}} \otimes \mathscr{O}_{X_{n-j}})
\end{equation*}
surjects for all $1\leq j\leq n-1$. Now we use the vector bundle technique (Lemma \hyperref[1.5]{1.5}) by taking $F=L\vert_{X_{n-j+1}}$, $R=L\vert_{X_{n-j+1}}$, $Q=\mathscr{O}_{X_{n-j+1}}(B\vert_{X_{n-j+1}})$, $r=1$, $G=B\vert_{X_{n-j}}$. We need to show the following:
\begin{equation}
    H^1(F\otimes Q^*)=0\textrm{; this comes from Lemma \hyperref[2.1]{2.1 (i)}.}\tag{3.3.5}
\end{equation}
\begin{equation}
    H^0(M_{L\vert_{X_{n-j}}}\otimes L\vert_{X_{n-j}})\otimes H^0(B\vert_{X_{n-j}})\rightarrow H^0(M_{L\vert_{X_{n-j}}}\otimes L\vert_{X_{n-j}}\otimes B\vert_{X_{n-j}} )\textrm{ surjects; it is  our hypothesis.}\tag{3.3.6}
\end{equation}
\begin{equation}
    H^0(L\vert_{X_{n-j}})\otimes H^0(B\vert_{X_{n-j}})\to H^0((L+B)\vert_{X_{n-j}})\textrm{ surjects; this comes from Lemma \hyperref[3.1]{3.1}.}\tag{3.3.7}
\end{equation}

That concludes the proof of the \textit{Induction Step}. Now we have to prove the base case.\par 
\noindent\textit{Base Case: We have to prove $H^0(M_{L\vert_{X_1}}\otimes L\vert_{X_1})\otimes H^0(B\vert_{X_1})\rightarrow H^0(M_{L\vert_{X_1}}\otimes L\vert_{X_1}\otimes B\vert_{X_1} )$ surjects.}\par
\noindent\textit{Proof of Base Case:} Notice, $deg(L\vert_{X_1})=(K+nB)\cdot B^{n-1}$, We have,
\begin{align*}
    2g-2=(B\vert_{X_2})^2+(B\vert_{X_2})\cdot K_{X_2}=B^n+(K+(n-2)B)\cdot B^{n-1}\textrm{ where $g=p_g(X_1)$.}
\end{align*}
\begin{align*}
    \implies deg(L\vert_{X_1})> 2g\textrm{, thanks to $B^n> 2$}\implies M_{L\vert_{X_1}}\textrm{ is semistable and $\mu(M_{L\vert_{X_1}})> -2$ (see \cite{Bu}).}
\end{align*}

We will use Proposition \hyperref[1.4]{1.4} to prove the required surjection of \textit{Base Case}.
We need to check:
\phantomsection\label{3.3.8}
\begin{equation*}
    \mu(M_{L\vert_{X_1}}\otimes L\vert_{X_1})> 2g.\tag{3.3.8}
\end{equation*}
\begin{equation*}
    \mu(M_{L\vert_{X_1}}\otimes L\vert_{X_1})> 4g-deg(B\vert_{X_1})-2h^1(B\vert_{X_1}).\phantomsection\label{3.3.9}\tag{3.3.9}
\end{equation*}

To prove \hyperref[3.3.8]{(3.3.8)}, we have to show $(K+nB)\cdot B^{n-1}-2\geq B^n+B^{n-1}\cdot (K+(n-2)B)+2$ which follows since $B^n\geq 4$.\par

Showing \hyperref[3.3.9]{(3.3.9)} is equivalent to proving $2h^1(B\vert_{X_1})\geq (n-3)B^n+B^{n-1}\cdot K+6$. Riemann-Roch gives the following, 
\begin{equation*}
    2h^1(B\vert_{X_1})=2h^0(B\vert_{X_1})+(n-3)B^n+B^{n-1}\cdot K.
\end{equation*}
That finishes the proof since $h^0(B\vert_{X_1})\geq 3$, thanks to $h^0(B)\geq n+2$.\QEDA\par

These three lemmas will help us prove the normal presentation of adjoint bundles associated to an ample, globally generated line bundle on a regular variety under suitable conditions.\par

\noindent\textit{{\bf Theorem 3.4.} \phantomsection\label{3.4}Let $X$ be a regular $n$-fold, $n\geq 3$. Let $B$ be an ample and base point free line bundle on $X$. We further assume:\\
\indent (a) $K$ is nef, $K+B$ is base point free.\\
\indent (b) $h^0(B)\geq n+2$.\\
\indent (c) $h^0(K+B)\geq h^0(K)+n+1$.\\
\indent (d) $(n-2)B-(n-1)K$ is nef and non-zero effective divisor.\\
Then $K+lB$ will satisfy the property $N_1$ for $l\geq n$.}\par

\noindent\textit{Proof.} We prove the assertion only for $l=n$, the case $l>n$ is similar. Let $L=K+nB$.
Since we already know that $H^1(M_L\otimes L)=0$ which comes from the projective normality of $L$, we only have to prove for all $k\geq 1$,\phantomsection\label{3.4.1} 
\begin{equation}
    H^1(M_{L}^{\otimes 2}\otimes L^{\otimes k})=0.\tag{3.4.1}
\end{equation}
\indent We omit the proof when $k\geq 2$ which follows easily from CM Lemma (Lemma \hyperref[1.7]{1.7}). Here we only prove the key case $k=1$ that is $H^1(M_{L}^{\otimes 2}\otimes L)=0$.
We have the short exact sequence:
\[
\begin{tikzcd}
 0 \arrow{r}{} & M_{L}^{\otimes 2}\otimes L \arrow{r}{}  & H^0(L)\otimes M_L\otimes L \arrow{r}{} & M_L\otimes L^{\otimes 2} \arrow{r}{} & 0.
\end{tikzcd}
\]
It is enough to prove that $H^0(L)\otimes H^0(M_L\otimes L)\rightarrow H^0(M_L\otimes L^{\otimes 2})$ surjects as $H^1(M_L\otimes L)=0$. We use Observation \hyperref[1.2]{1.2}; it is enough to prove the following:
\begin{equation}
    H^0(M_L\otimes L)\otimes H^0(B)\rightarrow H^0(M_L\otimes L\otimes B)\textrm{ surjects.}\phantomsection\label{3.4.2}\tag{3.4.2}
\end{equation}
\begin{equation}
    H^0(M_L\otimes L)\otimes H^0(lB)\rightarrow H^0(M_L\otimes L\otimes lB)\textrm{ surjects for all $l\geq 2$.}
    \phantomsection\label{3.4.3}\tag{3.4.3}
\end{equation}
\begin{equation}
    H^0(M_L\otimes (K+(2n-1)B)\otimes H^0(K+B)\rightarrow H^0(M_L\otimes (2K+2nB))\textrm{ surjects.}\phantomsection\label{3.4.4}\tag{3.4.4}
\end{equation}

We have proved \hyperref[3.4.2]{(3.4.2)} in Lemma \hyperref[3.3]{3.3}.\par

In order to prove \hyperref[3.4.3]{(3.4.3)}, we again use Observation \hyperref[1.2]{1.2}.
Therefore it is enough to prove that $H^0(M_L\otimes (K+lB))\otimes H^0(B)\rightarrow H^0(M_L\otimes (K+(l+1)B))$ surjects for $l>n$. To prove this our obvious choice is to use CM lemma (Lemma \hyperref[1.7]{1.7}).
First, we want to show that $H^1(M_L\otimes(K+(l-1)B))=0$ which is equivalent to showing the surjection of the following map:
\begin{equation}
    H^0(L)\otimes H^0(K+(l-1)B)\rightarrow H^0(L+K+(l-1)B).\phantomsection\label{3.4.5}\tag{3.4.5}
\end{equation}

If $l=n+1$ the this has already been proved in Theorem \hyperref[2.3]{2.3}, Step \hyperref[7]{7}.
If $l>n+1$ then in order to show the surjection of \hyperref[3.4.5]{(3.4.5)}, it is enough to prove $$H^0(2K+2nB+rB)\otimes H^0(B)\rightarrow H^0(2K+2nB+(r+1)B)$$ surjects for all $r\geq 0$.
This is \textit{Step} \hyperref[1]{1} in Theorem \hyperref[2.3]{2.3} with $k=2$.
Now we will show that, for all $2\leq i\leq n$, $H^i(M_L\otimes(K+(l-i)B))=0$.
We have the short exact sequence:
\[
\begin{tikzcd}
 0 \arrow{r}{} & M_{L}\otimes (K+(l-i)B) \arrow{r}{}  & H^0(L)\otimes (K+(l-i)B) \arrow{r}{} & 2K+(l+n-i)B \arrow{r}{} & 0.
\end{tikzcd}
\]
It gives us the long exact sequence: 
\[
\begin{tikzcd}[row sep=large, column sep=2ex]
 ... \arrow{r}{} & H^{i-1}(2K+(l+n-i)B) \arrow{r}{}  & H^i(M_{L}\otimes (K+(l-i)B)) \arrow{r}{} &  H^0(L)\otimes H^i(K+(l-i)B)\arrow{r}{} & ...
\end{tikzcd}
\]
Since the first and the last terms are zero by Kodaira vanishing, hence $H^i(M_L\otimes(K+(l-i)B))=0$ for all $2\leq i\leq n$.\par

We are left to prove \hyperref[3.4.4]{(3.4.4)}. Again we are going to use CM Lemma (Lemma \hyperref[1.7]{1.7}).
We have to prove the following three things:
\begin{equation}
    H^1(M_L\otimes(2n-2)B)=0.\phantomsection\label{3.4.6}\tag{3.4.6}
\end{equation}
\begin{equation}
    H^j(M_L\otimes ((2n-1-j)B-(j-1)K))=0\textrm{ for all $2\leq j\leq n-1$.}\phantomsection\label{3.4.7}\tag{3.4.7}
\end{equation}
\begin{equation}
    H^n(M_L\otimes ((n-1)B-(n-1)K))=0.\phantomsection\label{3.4.8}\tag{3.4.8}
\end{equation}

We observe that \hyperref[3.4.6]{(3.4.6)} is equivalent to showing $H^0(L)\otimes H^0((2n-2)B)\rightarrow H^0(L+(2n-2)B)$ surjects.
Using Observation \hyperref[1.2]{1.2}, this is equivalent to showing $H^0(K+lB)\otimes H^0(B)\rightarrow H^0(K+(l+1)B)$ surjects for all $l\geq n$.
This follows from Lemma \hyperref[2.1]{2.1} and Theorem \hyperref[2.3]{2.3}, \textit{Step} \hyperref[4]{4}.\par

To prove \hyperref[3.4.7]{(3.4.7)}, we write down the short exact sequence:\phantomsection\label{3.4.9}
\begin{equation}
\begin{tikzcd}[row sep=large, column sep=2ex]
 0 \arrow{r}{} & M_{L}\otimes (a_jB-b_jK) \arrow{r}{}  & H^0(L)\otimes (a_jB-b_jK) \arrow{r}{} & L\otimes (a_jB-b_jK) \arrow{r}{} & 0\tag{3.4.9}
\end{tikzcd}
\end{equation}
where $a_j=2n-1-j$, $b_j=j-1$.
The long exact sequence corresponding to it is:
\[
\begin{tikzcd}[row sep=large, column sep=2ex]
 ... \arrow{r}{} & H^{j-1}(L\otimes (a_jB-b_jK)) \arrow{r}{}  & H^j(M_{L}\otimes (a_jB-b_jK)) \arrow{r}{} & H^0(L)\otimes H^j(a_jB-b_jK) \arrow{r}{} & ...
\end{tikzcd}
\]
Now $H^{j-1}(L\otimes (a_jB-b_jK))=H^{j-1}(K+(3n-2j)B+(j-1)(B-K))=0$ as $3n-2j>0$ for all $j<n$ and $B-K$ is nef.
Also, $H^j(a_jB-b_jK)= H^j(K+(2n-2j-1)B+j(B-K))=0$ as $2n-2j-1>0$ for all $j<n$ and $B-K$ is nef.\par

We are left to prove \hyperref[3.4.8]{(3.4.8)} only.
The long exact sequence associated to \hyperref[3.4.9]{(3.4.9)} corresponding to $j=n$ gives the required vanishing for the following reasons:
\begin{equation*}
    H^{n-1}((2n-1)B-(n-2)K)=H^{n-1}(K+nB+(n-1)(B-K))=0.
\end{equation*}
\begin{equation*}
    H^n((n-1)(B-K))\myeq H^0(nK-(n-1)B)=H^0((n-1)K-(n-2)B+K-B)=0.
\end{equation*}
The last equality comes from the fact that $(n-1)K-(n-2)B$ is negative effective and $K-B$ is negative nef. That concludes the proof.\QEDA\par

Now we prove a weaker result for the normal presentation of the adjunction bundle associated to an ample, globally generated line bundle on an irregular variety of dimension $n$. Here we have to use Skoda complex to restrict ourselves to the multiplication map on the curve section as the variety is not regular. We include only a sketch of the proof as it is very similar to what we have done thus far.\par

\noindent\textit{{\bf Theorem 3.5.} \phantomsection\label{3.5}Let $X$ be an irregular variety of dimension $n$, $n\geq 3$. Let $B$ be an ample and base point free line bundle on $X$. We further assume:\\
\indent (a) $K$ is nef, $B'$ and $K+B'$ is base point free whenever $B\equiv B'$.\\
\indent (b) $h^0(B)\geq n+2$.\\
\indent (c) $h^0(K+B)\geq h^0(K')+n+1$ whenever $K\equiv K'$.\\
\indent (d) $(n-2)B-(n-1)K$ is nef and non-zero effective divisor.\\
Then $K+lB$ will satisfy the property $N_1$ for $l\geq n$.}\par

\noindent\textit{Proof.} As before, we just give the sketch for $L=K+nB$. We have to prove the vanishing $H^1(M_L^{\otimes 2}\otimes L^{\otimes k})=0$. Again, we just discuss the case when $k=1$.
It is enough to prove the surjection of the following map,
\begin{equation*}
    H^0(M_L\otimes L)\otimes H^0(L)\rightarrow H^0(M_L\otimes L^{\otimes 2}).
\end{equation*}

Let $E$ be a torsion line bundle in $Pic^0(X)$ which is not $n$ torsion. Such an $E$ exists as $Pic^0(X)$ is an abelian variety when $X$ is irregular. Note that $B+E$ is globally generated by assumption (a).
Observation \hyperref[1.2]{1.2} tells us it is enough to check the following three maps surject:
\begin{equation}
    H^0(M_L\otimes (K+nB))\otimes H^0(B+E)\rightarrow H^0(M_L\otimes ((n+1)B+E)).\phantomsection\label{3.5.1}\tag{3.5.1}
\end{equation}
\begin{equation}
    H^0(M_L\otimes (K+rB+E))\otimes H^0(B)\rightarrow H^0(M_L\otimes ((r+1)B+E))\textrm{ for $n+1\leq r\leq 2n-2$.}\phantomsection\label{3.5.2}\tag{3.5.2}
\end{equation}
\begin{equation}
    H^0(M_L\otimes (K+(2n-1)B+E))\otimes H^0(K+B-E)\rightarrow H^0(M_L\otimes (2K+2nB)).\phantomsection\label{3.5.3}\tag{3.5.3}
\end{equation}

To show \hyperref[3.5.1]{(3.5.1)} surjects, we use CM Lemma (Lemma \hyperref[1.7]{1.7}). We have to prove the following,
\begin{equation*}
    H^i(M_L\otimes (K+nB-iB-iE)=0\textrm{ for all $1\leq i\leq n$.} 
\end{equation*}

When $2\leq i\leq n-1$ this follows easily by multiplying the exact sequence \hyperref[0]{$(*)$} by suitable line bundle and then taking the cohomology.\par  
When $i=n$, same computation shows the vanishing once we see that $H^0(nE)=0$.\par 
To prove the vanishing for $i=1$, we need to show the surjection of the following map:
\begin{equation*}
    H^0(L)\otimes H^0(K+(n-1)B-E)\rightarrow H^0(2K+(2n-1)B-E).
\end{equation*}
By Observation \hyperref[1.2]{1.2}, Lemma \hyperref[2.1]{2.1} and Theorem \hyperref[2.3]{2.3}, \textit{Step} \hyperref[4]{4}, it is enough to prove the surjection of $H^0(K+(2n-2)B)\otimes H^0(K+B-E)\rightarrow H^0(2K+(2n-1)B-E)$.
Now, $K+B-E$ is base point free by our assumption. Let $C$ be a curve section of $K+B-E$.
Using Skoda complex (Definition \hyperref[1.8]{1.8}) and Lemma \hyperref[1.6]{1.6}, it is enough to check,
\begin{equation*}
    h^1(((2n-3)B+E)\vert_C)\leq h^0(K+B)-(n+1).
\end{equation*}
Now, $h^1(((2n-3)B+E)\vert_C)=h^0((nK-(n-2)B-nE)\vert_C)=h^0(nK-(n-2)B-nE)\leq h^0(K-nE)$ thanks to assumption (d). So, the inequality follows thanks to assumption (c).\par

\hyperref[3.5.2]{(3.5.2)} and \hyperref[3.5.3]{(3.5.3)} follows from CM Lemma (Lemma \hyperref[1.7]{1.7}) as well.\QEDA\par

\noindent\textbf{Remark 3.5.1.} We always have $h^0(K+B)\geq h^0(K')+n$ on any $n$-fold if $K+B$, $K'+B$, $B$ are ample, base point free and $K\equiv K'$.\par 
\noindent\textit{Proof.} We have $K'=K+\delta$ where $\delta$ is a numerically trivial line bundle.
By Riemann-Roch, $h^0(K+B)=h^0(K+B+\delta)$. The assertion follows from an argument similar to the proof of Remark \hyperref[2.2.01]{2.2.1}.\QEDB

\section{Properties \texorpdfstring{$N_0$}{} and \texorpdfstring{$N_1$}{} for pluricanonical series on three and four-folds}\label{KX}
In this section, we will concentrate on the behavior of pluricanonical series. First, we will prove a theorem whose corollaries will give us effective results on three and four folds. We again work on smooth varieties only.\par

\noindent\textit{{\bf Theorem 4.1.} \phantomsection\label{4.1} Let $X$ be an $n$ dimensional variety and let $B$ be an ample, globally generated line bundle on $X$. Let $L$ be a nef line bundle on $X$. Moreover, assume:\\
\indent (a) $(n-1)(B-L)-K$ is ample.\\
\indent (b) $B-K$ is ample.\\
\indent (c) $B+L$ is globally generated.\\
\indent (d) $h^0(K-L)\leq h^0(B)-(n+1)$.\\
Then $nB+L$ will be very ample and it will embed $X$ as a projectively normal variety.}\par

\noindent\textit{Proof.} Let $X_n$ be $X$ and let $X_{n-j}$ be a smooth irreducible $(n-j)$ fold chosen from the complete linear system of $|B\vert_{X_{n-j+1}}|$ by Bertini, for all $1\leq j\leq n-1$.
By adjunction, $K_{X_{n-j}}=(K+jB)\vert_{X_{n-j}}$ for all $0\leq j\leq n-1$.
We have to prove $H^0(k(nB+L))\otimes H^0(nB+L)\rightarrow H^0((k+1)(nB+L))$ surjects. Here we show the key case that is the case when $k=1$. The proof for $k\geq 2$ is similar.
We break the proof into a few steps.\par

\textit{Step 1: $H^0(nB+L)\otimes H^0(B)\rightarrow H^0((n+1)B+L)$ surjects.} 
We have the following diagram where $\mathscr{I}$ is the ideal sheaf of $X_1$ in $X$, $V$ is the cokernel of $H^0(B\otimes \mathscr{I})\rightarrow H^0(B)$:
\phantomsection\label{4.1.1}
\begin{equation}
\begin{tikzcd}[row sep=large, column sep=1.5ex]
 0\arrow[r] & H^0(nB+L)\otimes H^0(B\otimes \mathscr{I}) \arrow[r] \arrow[d] &  H^0(nB+L)\otimes H^0(B) \arrow[r] \arrow[d] &  H^0(nB+L)\otimes V \arrow[r]\arrow[d] & 0\\
0 \arrow[r] & H^0(((n+1)B+L)\otimes \mathscr{I}) \arrow[r] & H^0((n+1)B+L) \arrow[r] & H^0(((n+1)B+L)\vert_{X_1}) \arrow[r] & 0\tag{4.1.1}
\end{tikzcd}
\end{equation}
Note that $H^0((rB+L)\vert_{X_{n-j}})\rightarrow H^0((rB+L)\vert_{X_{n-j-1}})$ surjects for all $0\leq j\leq n-2$, $r\geq n$  because of the vanishing $H^1(((r-1)B+L)\vert_{X_{n-j}})=0$. Therefore the bottom horizontal sequence is exact. Note that the top row is exact as well.
The leftmost vertical map is surjective. Indeed, tensoring the following exact sequence (recall that $W$ is the span of $n-1$ general sections of $B$)
\begin{equation}
    0 \to \bigwedge\limits^{n-1}W\otimes B^{-(n-1)}\to \dots \to \bigwedge\limits^{2}W\otimes B^{-2} \to  W\otimes B^{-1} \to \mathscr{I} \to 0\phantomsection\label{4.1.2}\tag{4.1.2}
\end{equation}
by $(n+1)B+L$, we get the following exact sequence where $L'=(n+1)B+L$\[
\begin{tikzcd}
 0 \arrow{r}{f_{n-2}} & \bigwedge\limits^{n-1}W\otimes L'\otimes B^{-(n-1)}\arrow{r}{f_{n-1}} & \dots \arrow{r}{f_2} &  W\otimes L'\otimes B^{-1} \arrow{r}{f_1} & L'\otimes\mathscr{I} \arrow{r}{} & 0.
\end{tikzcd}
\]
Therefore, to see that the leftmost vertical map surjects, we just need $H^1(ker(f_1))=0$. The proof of this vanishing is similar to that of Lemma \hyperref[2.1]{2.1}. It comes from the following two facts :\par 

\indent \textit{Fact 1: $H^{j}(ker(f_j))=0\implies H^{j-1}(ker(f_{j-1}))=0$ for all $2\leq j\leq n-2$} (proof similar to Claim \hyperref[14]{1}, Lemma \hyperref[2.1]{2.1}).\\ 
\indent \textit{Fact 2: $H^{n-2}(L'-(n-1)B)= H^{n-2}(2B+L) = 0$} (using $B-K$ is ample and Kodaira vanishing).  \par

Therefore, to prove the assertion of this step, it is enough to show the surjection of the rightmost vertical map. We use Lemma \hyperref[1.6]{1.6} to prove that, we need the following inequality:
\begin{equation*}
    h^1((n-1)B+L)\vert_{X_1})\leq h^0(B)-(n+1).
\end{equation*}
But $h^1((n-1)B+L)\vert_{X_1})=h^0((K-L)\vert_{X_1})=h^0(K-L)$ (the last equality is due to the ampleness of $B+L-K$) which proves the assertion of this step because of our assumption (d).\par

\textit{Step 2: $H^0(rB+L)\otimes H^0(B)\rightarrow H^0((r+1)B+L)$ surjects for all $r\geq n+1$.}
This comes from CM Lemma (Lemma \hyperref[1.7]{1.7}).\par

\textit{Step 3: $H^0((2n-1)B+L)\otimes H^0(B+L)\rightarrow H^0(2nB+2L)$ surjects .}
This comes from CM Lemma (Lemma \hyperref[1.7]{1.7}) as well thanks to assumption (a).\QEDA\par

\noindent\textit{{\bf Corollary 4.2.} \phantomsection\label{4.2}Let $X$ be an $n$ dimensional variety, $n\geq 3$, with ample canonical bundle $K$. We further assume that $lK$ is globally generated for all $l\geq n+2$. Then the following will hold:\\
\indent (i) If $h^0(K)\leq h^0((n+2)K)-(n+1)$ then $n(n+2)K$ is very ample and it embeds $X$ as a projectively \indent\indent normal variety.\\
\indent (ii) If $h^0((n+2)K)>n+1$ then $(n(n+2)+1)K$ is very ample and it embeds $X$ as a projectively \indent\indent normal variety.\\
\indent (iii) $(n(n+2)+m)K$ is very ample and it embeds $X$ as a projectively normal variety for all $m\geq 2$.}\par

\noindent\textit{Proof of (i), (ii).} Follows directly from Theorem \hyperref[4.1]{4.1} with $B=(n+2)K$ and $L=0$, $K$ respectively.\par

\noindent\textit{Proof of (iii).} Let $s=n+2$. The proof is entirely based on CM Lemma (Lemma \hyperref[1.7]{1.7}). We give an outline here. We divide the proof into a few steps.\par

\textit{Step 1:} \textit{ $H^0((ns+m)K)\otimes H^0(sK)\rightarrow H^0(((n+1)s+m)K)$ surjects for all $m\geq 2$}.
This comes from CM Lemma (Lemma \hyperref[1.7]{1.7}).\par

\textit{Step 2:} \textit{ $H^0(((2n-1)s+m)K)\otimes H^0((s+m)K)\rightarrow H^0((2ns+2m)K)$ surjects for all $m\geq 2$}. To prove this, first notice that if $m\geq s$, then $m=as+b$ where $a\geq 1$ and $b<s$.
In that case, by Observation \hyperref[1.2]{1.2}, it is enough to show $H^0(((2ns+m+(a-1)s)K)\otimes H^0((s+b)K)\rightarrow H^0((2ns+2m)K)$ surjects for all $m\geq 2$ which comes from CM Lemma (Lemma \hyperref[1.7]{1.7}).
If $m<s$, we can directly use CM Lemma (Lemma \hyperref[1.7]{1.7}).\par

The above two steps shows the surjectivity of $H^0((ns+m)K)\otimes H^0((ns+m)K)\rightarrow H^0((2ns+2m)K)$. \\
Similar calculation shows $H^0(k(ns+m)K)\otimes H^0((ns+m)K)\rightarrow H^0((k+1)(ns+m)K)$ surjects for all $k\geq 2$. \QEDA\par





Now we combine our results with the base point freeness theorems on three and four folds (see \cite{EL} and \cite{Ka}). In particular, we will use Theorems \hyperref[1.9]{1.9}, \hyperref[1.11]{1.11} and  Corollary \hyperref[1.10]{1.10}. For the statement of the Riemann-Roch formula, we refer to \cite{Har}, Appendix A.\par

\noindent\textbf{Corollary 4.3.} \phantomsection\label{4.3}\textit{Let $X$ be a smooth projective three-fold with ample canonical bundle $K$. Then $nK$ is very ample and embeds $X$ as a projectively normal variety for all $n \geq 12$.} \par

\noindent\textit{Proof.} We have by Riemann-Roch that $\chi(D)+\chi(-D)=\displaystyle{\frac{-K\cdot D^2}{2}}+2\chi(\mathscr{O}_X)$. Hence, $K\cdot D^2$ is even for any divisor $D$. In particular $K^3$ is even.  By Corollary \hyperref[1.10]{1.10} we have that $4K$ is base point free.\par 
By CM Lemma (Lemma \hyperref[1.7]{1.7}), the corollary is obvious for $n \geq 14$ since $K$ is ample and we have Kodaira vanishing. For $n=13$ we use Theorem \hyperref[4.1]{4.1} with $L=K$. Conditions (a), (b), (c) are easily satisfied. We need to check that $h^0(4K) \geq 5$. We note that by Riemann-Roch (see the formula given in \cite{Har}, Appendix A, Exercise 6.7) and Remark \hyperref[1.12.1]{1.12.1} we have $h^0(4K) \geq 6+h^0(2K)$ and hence we are done.\par 
For the case $n=12$ we again apply Theorem \hyperref[4.1]{4.1} but now with $L=0$. Here we need to check the fact that $h^0(4K)\geq h^0(K)+4$. If $h^0(K)=0$ then we are done trivially since $4K$ is ample and base point free. If not then we know that $K$ is effective and hence $h^0(K) \leq h^0(2K)$. The required inequality comes from the inequality  $h^0(4K) \geq 6+h^0(2K)$.\QEDA\par

\noindent\textbf{Corollary 4.4.} \phantomsection\label{4.4}\textit{Let $X$ be a smooth projective three-fold with ample canonical bundle $K$. Then we have that the embedding by $nK$ for $n \geq 13$ is normally presented.} \par

\noindent\textit{Proof.} Suppose that $L=nK$. We note that the cases $n=3l+1$ with $l \geq 4$ normal presentation of $nK$ directly follows from Riemann-Roch and Theorem \hyperref[3.4]{3.4} for regular threefolds and \hyperref[3.5]{3.5} for irregular threefolds using $B=lK$ respectively. While using Theorem \hyperref[3.5]{3.5} we need to check the conditions. We only check the conditions (a) and (c) below. The other conditions follow directly from the Riemann-Roch formula (\cite{Har}, Appendix A, Exercise 6.7) once we note that $K\cdot c_2(X)\geq 0$.\par 

First we show that condition (a) holds. We have $B=lK$, $l \geq 4$.  Suppose $B'\equiv B$, then we have that $B'-K$ is ample and $(B'-K)^3 > 27$ (Using $K^3\geq 2$) and $(B'-K)\cdot C \geq 3$ and $(B'-K)^2\cdot S \geq 9$ for any curve $C$ and surface $S$ respectively. Hence $B'$ is base point free by Theorem \hyperref[1.9]{1.9}. Similar reasoning will show that $K+B'$ is base point free as well.\par 
Now show that condition (c) holds. Since $K+B$ is ample and base point free, if $h^0(K')=0$ we are done. Otherwise $K'$ is effective and $h^0(K') \leq h^0(2K')$. Note that $h^0(4K) \leq h^0((l+1)K)$. But since all higher cohomology of $2K'$ vanishes by Kodaira vanishing we have by Riemann Roch that $h^0(2K')$ depends only on the numerical class of $K'$ and hence $h^0(2K)=h^0(2K')$. So it is enough to show that $h^0(4K)-h^0(2K) \geq 4$ which we have shown in the proof of Corollary \hyperref[4.3]{4.3} (in fact $\geq 6$).\par

For other cases, it is enough to show that $H^1(M_L^{\otimes 2}\otimes L^{\otimes k})=0$ since we have already shown projective normality for $nK$ for $n \geq 13$. We only show the case $k=1$ since for $k \geq 2$ the proof follows from CM Lemma (Lemma \hyperref[1.7]{1.7}). We have the following exact sequence 
\[
\begin{tikzcd}
0 \arrow[r] & M_L^{\otimes 2}\otimes L \arrow[r] & H^0(L) \otimes M_L \otimes L \arrow[r] & M_L \otimes L^{\otimes 2} \arrow[r] & 0. 
\end{tikzcd}\] 
Taking cohomology we have the following 
\[
\begin{tikzcd}
... \arrow[r] & H^0(L) \otimes H^0(M_L \otimes L) \arrow[r] & H^0(M_L \otimes L^{\otimes 2}) \arrow[r] & H^1(M_L^{\otimes 2}\otimes L) \arrow[r] & ...
\end{tikzcd}\]
It is enough to show that $H^0(L) \otimes H^0(M_L \otimes L) \rightarrow H^0(M_L \otimes L^{\otimes 2})$ is surjective. Now $4K$ is base point free. We first show that $H^0(M_L \otimes L) \otimes H^0(4K) \rightarrow H^0(M_L \otimes L+4K)$ is surjective. To do this it is enough to show (by Lemma \hyperref[1.7]{1.7}) the following three vanishings: 
\begin{equation}
    H^1(M_L \otimes L-4K) = 0,\tag{i}
\end{equation}
\begin{equation}
    H^2(M_L \otimes L-8K) = 0,\tag{ii}
\end{equation}
\begin{equation}
    H^3(M_L \otimes L-12K) = 0.\tag{iii}
\end{equation}

Now $L=nK$ with $n \geq 14$ (the case when $L=13K$ has already been taken care of). By tensoring the exact sequence
\[
\begin{tikzcd}
0 \arrow[r] & M_L \arrow[r] & H^0(L) \otimes \mathscr{O}_X \arrow[r] & L \arrow[r] & 0 
\end{tikzcd}\] 
by $L-8K$ and $L-12K$ respectively and using Kodaira vanishing theorem we can see that (ii) and (iii) follow immediately. Now we note that to show (i) we need to to show that the following map
\begin{equation*}
    H^0(L) \otimes H^0(L-4K) \rightarrow H^0(2L-4K)
\end{equation*}
is surjective. We now note that $L-4K = mK$ where $m \geq 10$. Using observation \hyperref[1.2]{1.2} we can keep on showing surjectivity of multiplication maps by $H^0(4K)$ until we are left with $lK$ where $0 \leq l \leq 7$. Hence $H^0(L) \otimes H^0(L-4K) \rightarrow H^0(2L-4K)$ is surjective for $L=nK$ and $n \geq 19$. We need to check separately from $14 \leq n \leq 18$.\par

\textit{Case n=14.} We need to show the surjectivity of $H^0(14K) \otimes H^0(10K) \rightarrow H^0(24K)$. 
We have by Lemma \hyperref[1.7]{1.7} the surjectivity of $H^0(14K) \otimes H^0(4K) \rightarrow H^0(18K)$. 
We have the surjectivity of $H^0(18K) \otimes H^0(6K) \rightarrow H^0(24K)$ using \textit{Step 1}, Theorem \hyperref[4.1]{4.1} with $B=6K$ and $L=0$. \par

\textit{Case n=15.} We need to show the surjectivity of $H^0(15K) \otimes H^0(11K) \rightarrow H^0(26K)$.
We have that $H^0(15K) \otimes H^0(5K) \rightarrow H^0(20K)$ surjects by Theorem \hyperref[4.1]{4.1} with $B=5K$ and $L=0$. We also have the surjectivity of $H^0(20K) \otimes H^0(6K) \rightarrow H^0(26K)$ by Lemma \hyperref[1.7]{1.7}. \par
The case $n=16$ is obvious.\par 
\textit{Case n=17.} We need to show the surjectivity of $H^0(17K) \otimes H^0(13K) \rightarrow H^0(30K)$.
This case is easy and follows from Lemma \hyperref[1.7]{1.7}. \par

\textit{Case n=18.} We need to show the surjectivity of $H^0(18K) \otimes H^0(14K) \rightarrow H^0(32K)$.
This case follows directly from Lemma \hyperref[1.7]{1.7}. \par

The algorithmic nature of the proof shows that we have actually proved the surjectivity of the map $H^0(M_L \otimes (L+4lK)) \otimes H^0(4K) \rightarrow H^0(M_L \otimes (L+4(l+1)K))$.
Since $L=nK$ where $n \geq 14$, to complete the proof we just need to prove the surjection of the multiplication map
\begin{equation*}
    H^0(M_L \otimes (L+4lK)) \otimes H^0(pK) \rightarrow H^0(M_L \otimes (L+(4l+p)K))
\end{equation*}
where $l \geq 2$ and $p \leq 7$. Moreover if $n\geq 16$ we have that $l \geq 3$. So for $n \geq 16$, using Lemma \hyperref[1.7]{1.7} we see that it is enough to prove the surjection of $H^0(L+mK) \otimes H^0(L) \rightarrow H^0(2L+mK)$
where $m\geq 5$. But we have the surjection of $H^0(L) \otimes H^0(L) \rightarrow H^0(2L)$. Thus, using Observation \hyperref[1.2]{1.2}, we only need to prove the surjection of $H^0(lK) \otimes H^0(mK) \rightarrow H^0((m+l)K)$ where $l\geq 32$.
Since $4K$ is base point free, we have the above surjection by Lemma \hyperref[1.7]{1.7} and Observation \hyperref[1.2]{1.2}.\par 
To finish the proof we need to handle the two following cases separately. \par

\textit{L=14K:} We need to show the surjection of 
$H^0(M_L \otimes (L+8K)) \otimes H^0(6K) \rightarrow H^0(M_L \otimes (L+14K))$.\\
By lemma \hyperref[1.7]{1.7} we notice that it is enough to show the surjection of $H^0(16K) \otimes H^0(14K) \rightarrow H^0(30K)$
which is clear by the same lemma. \par

\textit{L=15K:} We need to show the surjection of 
$H^0(M_L \otimes (L+8K)) \otimes H^0(7K) \rightarrow H^0(M_L \otimes (L+15K))$.
By Lemma \hyperref[1.7]{1.7}, we notice that it is enough to show that $H^0(16K) \otimes H^0(15K) \rightarrow H^0(31K)$ surjects which is again clear by the same lemma.\QEDA\par

\noindent\textit{{\bf Corollary 4.5.} \phantomsection\label{4.5}Let $X$ be a smooth projective four dimensional variety with ample canonical bundle $K$. Then $nK$ is very ample and it will embed $X$ as a projectively normal variety for all $n\geq 24$.}\par

\noindent\textit{Proof.} It comes from Corollary \hyperref[4.2]{4.2}. 
The following is the Riemann-Roch formula a line bundle $B$,
\begin{equation*}
    \chi(B)=-\frac{1}{720}(K^4-4K^2\cdot c_2-3c_2^2+K\cdot c_3+c_4)-\frac{1}{24}B\cdot K\cdot c_2+\frac{1}{24}B^2\cdot (K^2+c_2)-\frac{1}{12}B^3\cdot K+\frac{1}{24}B^4.
\end{equation*}
It is enough to show that $h^0(2K)\leq h^0(6K)-5$ which can be seen easily, thanks to the fact $K^2\cdot c_2\geq 0$ (see Remark \hyperref[1.12.1]{1.12.1}). In fact, $h^0(2K)\leq h^0(6K)-6$ which verifies condition (ii).\QEDA\par

\noindent\textbf{Corollary 4.6.} \phantomsection\label{4.6}\textit{Let $X$ be a smooth projective 4-fold with ample canonical bundle $K$ we have that the embedding by $nK$ for $n \geq 25$ is normally presented.} \par

\noindent\textit{Proof.} We use the same argument as in Corollary \hyperref[4.5]{4.5}, but now using the fact that $6K$ is globally generated (see Theorem \hyperref[1.11]{1.11}).\QEDA\par
 
\noindent\textbf{Corollary 4.7.} \phantomsection\label{4.7}\textit{Let $X$ be a smooth projective 5-fold with ample canonical bundle $K$ with an additional property that $p_g(X)\geq 1$. Then the embedding by $nK$ for $n \geq 35$ is projectively normal and the embedding by $nK$ for $n\geq 36$ is normally presented.} \par

 \noindent\textit{Proof.} We know that $nK$ is globally generated for $n\geq 7$ (see \cite{Zhu}). Let $\mathscr{K}$ be a smooth divisor chosen from the linear system of $|K|$. The corollary follows at once if we notice the fact that $h^0(7K)-h^0(6K)=h^0((7K)\vert_\mathscr{K})$ and apply Riemann-Roch formula on $\mathscr{K}$ to verify conditions (i) and (ii) of Corollary \hyperref[4.2]{4.2}. Normal Presentation follows from the similar arguments used before.\QEDA

\section{Appendix}\label{AP}
We list a remark and its corollaries. The remark discusses the case when the variety is singular. 

 
 \noindent\textbf{Remark 5.1.} \phantomsection\label{6.2}We note that the fact that the Theorems \hyperref[2.3]{2.3}, \hyperref[3.4]{3.4}, \hyperref[3.5]{3.5} and \hyperref[4.1]{4.1} goes through for $X$ normal, Cohen-Macaulay with Du Bois singularities. We discuss below the precise reasons for which we require smoothness.\par 
  We require smooth hyperplane sections of the ample and base point free line bundle $B$. The smoothness is used to justify Kodaira vanishing (both on the general member of $|B|$ and on $X$) and to apply Green's result (Lemma \hyperref[1.6]{1.6}) on the smooth curve section. \par 
 We observe that if $X$ is normal, Cohen-Macaulay with Du Bois singularities and $B$ is Cartier, the general member of $|B|$ is Cohen-Macaulay as well. Also since $X$ is nonsingular in codimension $1$ and $|B|$ is base point free, a general member of $|B|$ is smooth outside the singular locus of $X$ (by Bertini's theorem) and is hence nonsingular in codimension $1$. The above two observations show that the general member is normal. Now the general hyperplane section of $|B|$ also has Du Bois singularities (see \cite{Kol}, Proposition 6.20). We also have (see \cite{Kol}, Theorem 10.42) that for a projective, Cohen-Macaulay variety with Du Bois singularities, Kodaira vanishing theorem holds for an ample line bundle.
Now the complete intersection surface that we get is a normal surface and hence singularities are isolated. So Bertini's Theorem gives us a smooth curve section and we can apply Lemma \hyperref[1.6]{1.6}. \QEDB\par

Now Koll\'ar and Kov\'acs prove that log canonical singularities are Du Bois (see \cite{Kol1}) and hence by Remark \hyperref[6.2]{5.1}, we have that the results mentioned in Remark \hyperref[6.2]{5.1} go through for  log canonical singularities and hence in particular for $\mathbb{Q}$-factorial terminal Gorenstein and canonical Gorenstein singularities. Oguiso-Peternell's generalization of Ein-Lazarsfeld's result on Fujita conjecture combined with Theorem \hyperref[3.4]{3.4}, Theorem \hyperref[3.5]{3.5} and Theorem \hyperref[4.1]{4.1} gives the following.\par

\noindent\textbf{Corollary 5.2.} \textit{If $X$ is a projective three-fold with $\mathbb{Q}$-factorial terminal Gorenstein singularities and ample canonical bundle $K$ then we have that $18K$ is projectively normal and $19K$ is normally presented.} \par
\noindent\textit{Proof.} By \cite{OP}, Corollary 2.3, we have that $6K$ is base point free and hence the result follows by exact same argument as in Corollary \hyperref[4.3]{4.3} of this article provided we have the exact form of Riemann-Roch as on a smooth threefold. Now, by Theorem 10.2, \cite{Reid}, we have what we want since the Gorenstein assumption gives us that $6K$ is Cartier and hence we do not have any contribution due to the singularities of the sheaf $\mathscr{O}_X(6K)$.\QEDA \par

\noindent\textbf{Corollary 5.3.} \textit{If $X$ is a projective three-fold with canonical Gorenstein singularities and ample canonical bundle $K$ then we have that $24K$ is projectively normal and $25K$ is normally presented.} \par

\noindent\textit{Proof.} The result follows since $8K$ is base point free (see \cite{OP}, Corollary 2.2).\QEDA
 
\bibliographystyle{plain}

\end{document}